\documentclass[reqno,11pt]{amsart}
\usepackage{psfig, amsmath, amsfonts, amssymb, amsthm, amscd}

\DeclareMathSymbol{\leqslant}{\mathalpha}{AMSa}{"36} 
\DeclareMathSymbol{\geqslant}{\mathalpha}{AMSa}{"3E} 
\DeclareMathSymbol{\eset}{\mathalpha}{AMSb}{"3F}     
\renewcommand{\leq}{\;\leqslant\;}                   
\renewcommand{\geq}{\;\geqslant\;}                   
\DeclareMathOperator*{\inter}{\bigcap}       
\newcommand{\limtwo}[2]{\lim_{\substack{#1 \\ #2}}}     

\makeatletter
\def\captionfont@{\footnotesize}
\def\captionheadfont@{\scshape}

\long\def\@makecaption#1#2{%
  \vspace{2mm}
  \setbox\@tempboxa\vbox{\color@setgroup
    \advance\hsize-6pc\noindent
    \captionfont@\captionheadfont@#1\@xp\@ifnotempty\@xp
        {\@cdr#2\@nil}{.\captionfont@\upshape\enspace#2}%
    \unskip\kern-6pc\par
    \global\setbox\@ne\lastbox\color@endgroup}%
  \ifhbox\@ne 
    \setbox\@ne\hbox{\unhbox\@ne\unskip\unskip\unpenalty\unkern}%
  \fi
  \ifdim\wd\@tempboxa=\z@ 
    \setbox\@ne\hbox to\columnwidth{\hss\kern-6pc\box\@ne\hss}%
  \else 
    \setbox\@ne\vbox{\unvbox\@tempboxa\parskip\z@skip
        \noindent\unhbox\@ne\advance\hsize-6pc\par}%
\fi
  \ifnum\@tempcnta<64 
    \addvspace\abovecaptionskip
    \moveright 3pc\box\@ne
  \else 
    \moveright 3pc\box\@ne
    \nobreak
    \vskip\belowcaptionskip
  \fi
\relax
}
\makeatother
\def\writefig#1 #2 #3 {\rlap{\kern #1 truecm
\raise #2 truecm \hbox{#3}}}
\def\figtext#1{\smash{\hbox{#1}}
\vspace{-5mm}}
\psfigurepath{.:./pictures}

\newtheorem{lem}{Lemma}[section]
\newtheorem{pro}{Proposition}[section]
\newtheorem{thm}{Theorem}[section]

\newtheorem{remark}{Remark}[section]

\newcommand{\cA}{\ensuremath{\mathcal A}}
\newcommand{\cB}{\ensuremath{\mathcal B}}
\newcommand{\cC}{\ensuremath{\mathcal C}}

\newcommand{\cE}{\ensuremath{\mathcal E}}

\newcommand{\cH}{\ensuremath{\mathcal H}}

\newcommand{\cK}{\ensuremath{\mathcal K}}
\newcommand{\cL}{\ensuremath{\mathcal L}}
\newcommand{\cM}{\ensuremath{\mathcal M}}

\newcommand{\cO}{\ensuremath{\mathcal O}}
\newcommand{\cP}{\ensuremath{\mathcal P}}

\newcommand{\cV}{\ensuremath{\mathcal V}}
\newcommand{\cW}{\ensuremath{\mathcal W}}
\newcommand{\cX}{\ensuremath{\mathcal X}}

\newcommand{\frB}{\ensuremath{\mathfrak B}}

\newcommand{\frI}{\ensuremath{\mathfrak I}}
\newcommand{\frJ}{\ensuremath{\mathfrak J}}

\newcommand{\frW}{\ensuremath{\mathfrak W}}

\newcommand{\frg}{\ensuremath{\mathfrak g}}

\newcommand{\bbB}{{\ensuremath{\mathbb B}} }

\newcommand{\bbD}{{\ensuremath{\mathbb D}} }

\newcommand{\bbL}{{\ensuremath{\mathbb L}} }

\newcommand{\bbN}{{\ensuremath{\mathbb N}} }

\newcommand{\bbP}{{\ensuremath{\mathbb P}} }

\newcommand{\bbR}{{\ensuremath{\mathbb R}} }
\newcommand{\bbS}{{\ensuremath{\mathbb S}} }

\newcommand{\bbZ}{{\ensuremath{\mathbb Z}} }

\newcommand{\ga}{\alpha}
\newcommand{\gb}{\beta}
\newcommand{\gd}{\delta}
\newcommand{\gep}{\varepsilon}       

\newcommand{\gr}{\rho}
\newcommand{\gz}{\zeta}
\newcommand{\gG}{\Gamma}

\newcommand{\gD}{\Delta}

\newcommand{\go}{\omega}

\newcommand{\gs}{\sigma}

\newcommand{\Ham}{{\bf H}}         
\newcommand{\PF}{{\bf Z}}          

\newcommand{\Is}{\mu}              
\newcommand{\Perc}{\Phi}           
\newcommand{\Joint}{\bbP}          

\newcommand{\Boxx}[1]{\bbD_{#1}}   
\newcommand{\uBox}{\widehat\bbD}         

\newcommand{\normI}[1]{\|#1\|_{{\scriptscriptstyle 1}}}

\newcommand{\bdf}{\eta}

\newcommand{\intSTbdf}[1]{\cW_{\gb,#1}}

\newcommand{\smallo}{o}

\newcommand{\bnd}{\partial}
\newcommand{\abs}[1]{\lvert#1\rvert}

\newcommand{\ra}{\rightarrow}
\newcommand{\setof}[2]{\{#1\,:\,#2\}}

\newcommand{\BoxxN}{{\Boxx{N}}}
\newcommand{\pairs}[2]{(#1\,,\,#2)}

\newcommand{\halfspace}{{\bbL^d}}

\newcommand{\taubd}{\Delta_{\gb,\eta}}
\newcommand{\Isbd}[2]{\Is^{\beta,\bdf}_{#2,#1}}

\newcommand{\Isbdf}[3]{\Is^{\beta,#3}_{#2,#1}}
\newcommand{\PFbd}[2]{\PF^{\beta,\bdf}_{#2,#1}}

\newcommand{\vol}{\text{vol}}

\newcommand{\RCbdf}[2]{\Perc^{\gb,#2,#1}_N}

\newcommand{\RCbdpmf}[1]{\Perc^{\gb,#1}_{N,\pm}}
\newcommand{\wired}{{\scriptscriptstyle\rm w}}
\newcommand{\jp}{\Joint_{N}^{\gb,\bdf}}
\newcommand{\jpm}{\Joint_{N,\pm}^{\gb,\abs\bdf}}

\def\1{\ifmmode {1\hskip -3pt \rm{I}} \else {\hbox {$1\hskip -3pt \rm{I}$}}\fi}
\newcommand{\BV}{{ {\rm BV}({\rm int}\uBox,\{\pm1\}) }}         
\def\lra{\leftrightarrow}
\def\nlra{\not \leftrightarrow}
\def\vol{{\rm vol}}
\newcommand{\df}{\stackrel{\gD}{=}}
\newcommand{\sDor}[1]{\widehat{\bbD}_{#1}}   
\newcommand{\sBox}[1]{\widehat{\bbB}_{#1}}      
\newcommand{\dBox}[1]{{\bbB}_{#1}} 

\newcommand{\st}{\tau_{\gb}}
\newcommand{\lb}{\left(}
\newcommand{\rb}{\right)}

\newcommand{\IsN}{\Is_{N,+}^{\gb}}
\newcommand{\FKm}[2]{\Phi^{#1}_{#2}}
\newcommand{\FKpf}[2]{\PF^{#1}_{#2}}

\newcommand{\symdiff}{{\scriptstyle \triangle}\,}

\begin{document}
\title[\,]{Winterbottom construction for finite range ferromagnetic
  models: An ${\Bbb L}_1$-approach}
\author{T.~Bodineau}
\address{
Universit\'e Paris 7, D\'epartement de Math\'ematiques, Case 7012, 2 place
Jussieu, F-75251 Paris, France}
\email{Thierry.Bodineau\@@gauss.math.jussieu.fr}
\author{D.~Ioffe}
\address{
Faculty of Industrial Engineering, Technion, Haifa 32000, Israel}
\email{ieioffe\@@ie.technion.ac.il}
\thanks{Partly supported by the ISRAEL SCIENCE FOUNDATION founded by The Israel
Academy of Science and Humanities and by the
 Technion V.P.R fund GLASBERG-KLEIN research fund}
\author{Y.~Velenik}
\address{
Laboratoire d'Analyse, Topologie et Probabilit\'es, UMR-CNRS 6632, CMI,
Universit\'e de Provence, 39 rue Joliot Curie, 13453 Marseille, France}
\email{velenik\@@cmi.univ-mrs.fr {\sl Homepage:}
http://www.cmi.univ-mrs.fr/$\scriptstyle\sim$velenik/}
\thanks{
Partly supported by a Swiss National Science Foundation grant
\#8220-56599. 
}
\date{\today}
\begin{abstract}
We provide a rigorous microscopic derivation of the thermodynamic description
of equilibrium crystal shapes in the presence of a substrate, first studied by
Winterbottom.
We consider finite range ferromagnetic Ising models with pair interactions in
dimensions greater or equal to $3$, and model the substrate by a finite-range
boundary magnetic 
 field acting on the spins close to the bottom wall of the box. 
\end{abstract}
\maketitle


\vspace*{5mm}

\section{Introduction and Results}
\subsection{Winterbottom variational problem}
\label{subsection_intro}
The Winterbottom theory \cite{Wi}  gives a phenomenological prediction for the 
equilibrated shapes of small (that is disregarding gravitation) crystal
particles placed on solid substrates. The notion of equilibrium comes on
two levels: 

\noindent
(i) As far as the bulk properties are considered, both the particle and the
vapor around it are assumed to  represent two different {\em thermodynamic
 phases} of the same physical substance.
 
\noindent   
(ii) The total interfacial free energy of the system comprises the  anisotropic
{\em surface tension} between the particle and the vapor, as well as the  {\em
boundary surface tensions} between the particle and the  substrate and,
respectively,  between the vapor and the  substrate. The substrate could,
thus, exhibit  a preference towards one of the  phases, and the macroscopic 
Winterbottom equilibrium shape corresponds to the state of  minimal 
interfacial energy at a given particle volume.
\begin{figure}[t]
\centerline{\raise .5 Truecm \hbox{\psfig{file=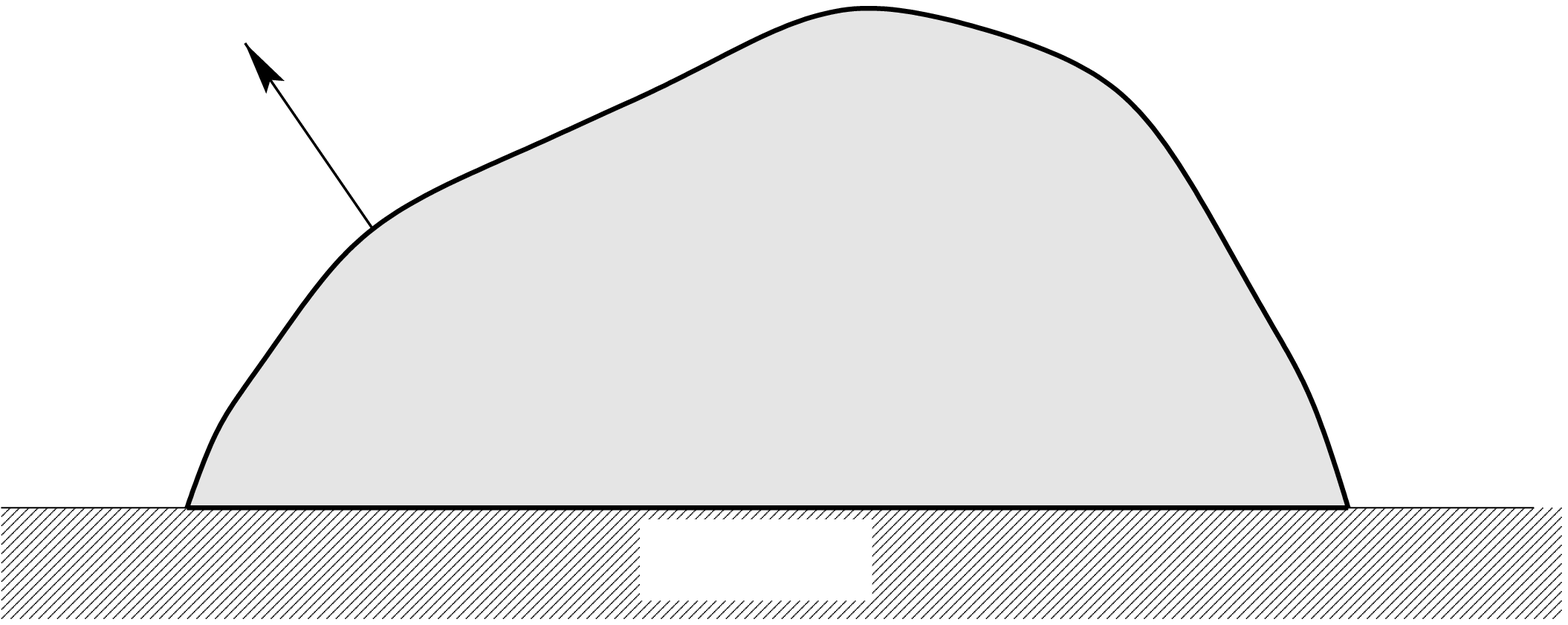,height=2cm}
}
\hspace*{2cm}
\psfig{file=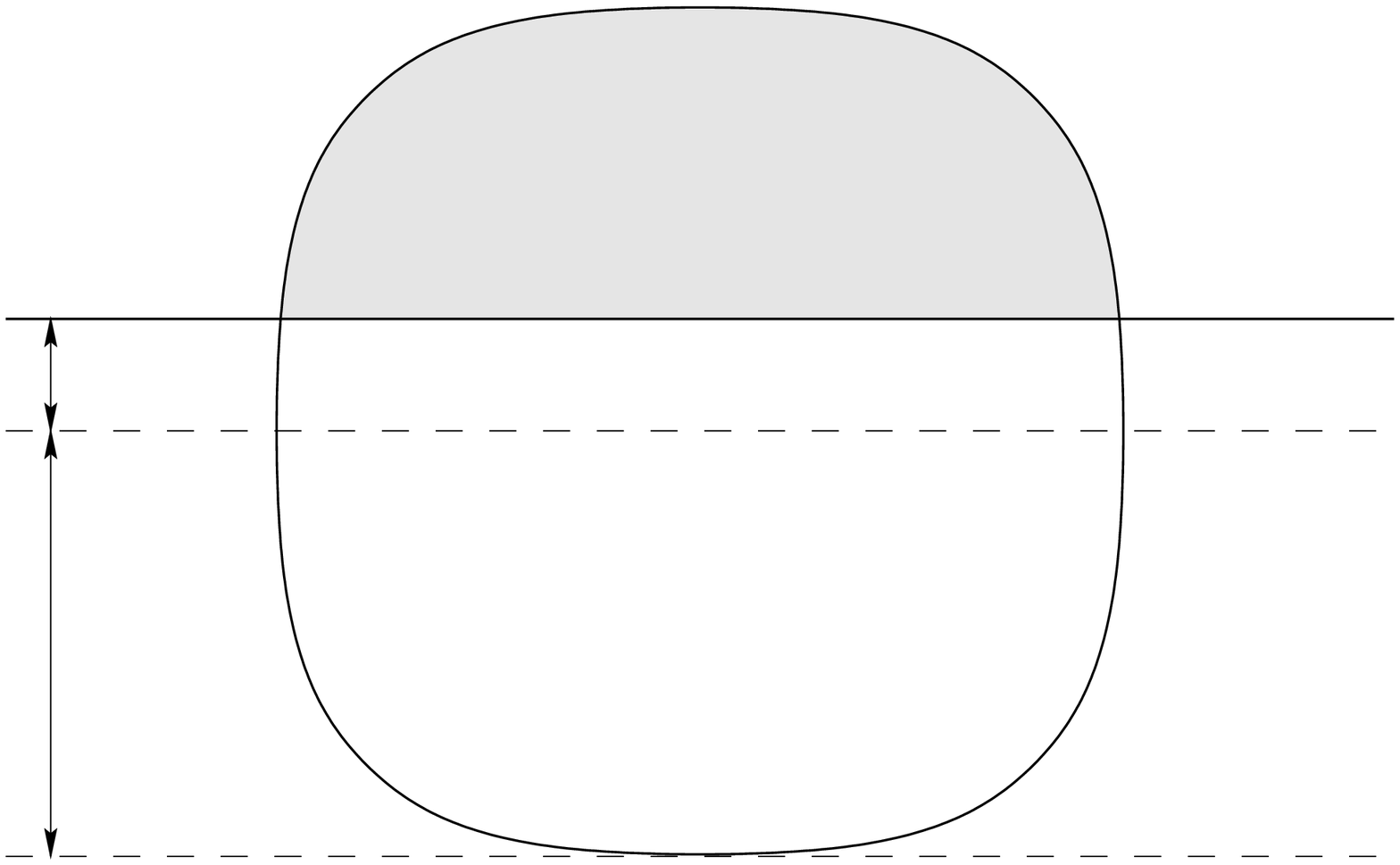,height=3cm}}
\hspace*{-0.25cm}
\figtext{ 
\writefig       -5.50   2.35    {\footnotesize $\vec{n}_x$} 
\writefig       -4.85   2.00    {\footnotesize $x$} 
\writefig       -3.50   2.00    {\footnotesize $P$}
\writefig       -3.95   1.00    {\footnotesize $\partial P_{\rm SC}$}
\writefig       -3.00   3.00    {\footnotesize $\partial P_{\rm CV}$} 
\writefig       -2.00   3.50    {\footnotesize Vapour} 
\writefig        0.50   2.10    {\footnotesize $-\gD_{\gb,\bdf}$} 
\writefig        0.50   1.15    {\footnotesize $\tau_\gb(\vec{e}_d)$} 
}
\caption{Left: The crystal $P$ on the substrate. The boundary is decomposed into two
pieces, corresponding to the crystal/vapor interface and the crystal/substrate
interface. Right: The Winterbottom shape $\cK_{\beta ,\bdf}$ 
(shaded region) is obtained by intersecting the Wulff shape and the half-plane; on
the picture $\gD_{\gb ,\bdf}$ is negative.}
\label{fig_energy}
\end{figure} 

\vskip 0.1in
We use $\beta$ to denote the inverse temperature at which the two phases
(crystal and vapor) coexist, and the parameter $\bdf$ to quantify the 
preference of the substrate towards one of the two pure phases. In the 
 magnetic interpretation of the microscopic models we shall consider below $\bdf$ 
is the boundary magnetic field. 
  According to the suppositions (i),(ii) above, the
two principal players to determine the equilibrium shape  are the  anisotropic
crystal-vapor surface tension $\tau_\beta =\tau_{\beta
,\text{CV}}:{\bbS}^{d-1}\mapsto {\bbR}_+$ and  the  number $\Delta_{\beta ,\bdf}
=\tau_{\beta ,\text{SC}}-\tau_{\beta ,\text{SV}}$, where  $\tau_{\beta
,\text{SC}}$ and $\tau_{\beta ,\text{SV}}$ are, respectively, the 
substrate-crystal and the substrate-vapor boundary surface tensions. Indeed,
if the crystal particle occupies the region $P$ of ${\bbR}^d$  its boundary
splits into two disjoint pieces  $\partial P=\partial
P_{\text{SC}}\bigcup\partial P_{\text{CV}}$ (Figure~\ref{fig_energy}), and, up
to a constant, the total surface  energy of the system can be written (at least
for piece-wise smooth shapes) as
\begin{equation}
\label{intro_functional}
{\cW}_{\beta ,\bdf}\left( P\right)~\df~\int_{\partial P_{\text{CV}}}\tau_\beta (\vec{n}_{x})
\text{d}S(x) ~+~\Delta_{\beta ,\bdf} \, \left| \partial P_{\text{SC}}
\right| .
\end{equation}
The equilibrium particle shapes at a prescribed volume $v$ correspond, in  the
above notation, to the  minimizers of the following variational problem:\\[2mm]
$\text{{\bf (WP)}}_{\beta ,\bdf}\qquad\qquad\qquad {\cW}_{\beta ,\bdf}\left(
P\right)~\longrightarrow~\text{min} \qquad\text{Given}\ \ \vol_d\left(
P\right)=v .$\\[2mm]
%
Winterbottom \cite{Wi} described the geometric construction of the solutions to
$\text{\bf (WP)}_{\beta ,\bdf}$, which, in fact, is an easy consequence of the
general Wulff construction \cite{Wu} of the equilibrium crystal shapes: First 
of all the minimizers are scale invariant under dilatation. In order to
construct the unnormalized Winterbottom shape pick the ``free'' Wulff crystal
$$
\cK_{\beta}~\df~
\bigcap_{n\in{\bbS^{d-1}}}\left\{ x\in \bbR^d ~\Big|~(x,n)\leq \tau_\beta (n)\right\} \, ,
$$
and intersect it with the half-space
\begin{equation} 
\label{0_shape}
\cK_{\beta ,\bdf}~\df~\cK_{\beta}\bigcap\left\{
x\in \bbR^d ~\Big|~x_d\geq -\Delta_{\beta ,\bdf}\right\},
\end{equation}
where we use the coordinate representation $x =(x_1 ,...,x_d )$ of a vector
$x\in\bbR^d$. It is convenient to parameterize the substrate surface by one of
the  coordinate hyperplanes, say by $\cP_d\df \{ x : x_d =0\}$.  The
equilibrium  Winterbottom shape at volume $v$ is, then, given by
\begin{equation} 
\label{1_shape}
\cK_{\beta ,\bdf}(v)~\df~
\left(\frac{v}{\lambda}\right)^{1/d}\left( \Delta_{\gb,\bdf}\vec{e}_d +
\cK_{\beta ,\bdf}\right) ,
\end{equation}
where $\vec{e}_d$ is the unit coordinate vector in the $d$-direction and 
$\lambda =\lambda (\beta ,\bdf )=\vol_d\left(\cK_{\beta ,\bdf}\right)$ is used
to denote  the volume of the unnormalized Winterbottom shape.

As we shall explain in Subsection~\ref{subsection_wallfree energy},
 the value of the difference
$\Delta_{\beta ,\bdf}$ always lies in the  interval $\left[-\tau_\beta (\vec{e}_d ),
\tau_\beta (\vec{e}_d )\right]$\footnote{We shall always assume that
$\st(x)=\st(-x)$.}.  Thus,  one of the following three cases
happens:
\begin{enumerate}
\item {\em Complete drying:} $\Delta_{\beta ,\bdf} = \tau_\beta (\vec{e}_d)$. In this
case   $\cK_{\beta ,\eta} = \cK_\beta$, which means that the
Winterbottom droplet has no energetic reasons to be near the substrate surface 
and can  be located anywhere above it.
\item {\em Partial drying/wetting:} In this regime 
$\Delta_{\beta ,\bdf}\in \left(-\tau_\beta (\vec{e}_d) ,\tau_\beta
(\vec{e}_d\right))$. Accordingly, the convex body $\cK_{\beta ,\bdf}$ is
obtained by removing a proper cap from  $\cK_\beta$  (Figure
\ref{fig_energy}). We shall, furthermore, distinguish between the partial 
drying, $\Delta_{\beta ,\bdf}>0$, and the partial wetting, $\Delta_{\beta ,\bdf}\leq 0$,
cases.  
\item {\em Complete wetting:} $\Delta_{\gb ,\bdf}= -\tau_\beta (\vec{e}_d)$. In this
case, $\cK_{\beta ,\bdf}$ has volume  zero and, therefore, cannot be rescaled
to get a volume $v$.  Physically this corresponds to the formation of a
microscopic film along the substrate,  preventing the vapor phase to reach it.
\end{enumerate}

Our main objective here is to give a rigorous derivation  of  the
phenomenological Winterbottom picture 
 in the scaling limit  of microscopic models of
lattice gases, that is directly from the basic principles of the Equilibrium
Statistical Mechanics. The class of microscopic models we consider are
ferromagnetic Ising type models with finite range of interactions. In the
subsequent subsections we shall briefly recall these models as well as the
thermodynamic procedures  leading to the notion of {\em bulk phases} and
various {\em surface tensions} which appear in the functional
\eqref{intro_functional}. The results fall in the framework of the
${\bbL}_1$-theory of phase segregation, which has been  initiated in the works
\cite{Pisztora} on the nearest neighbor Ising model and in the  works
\cite{ABCP}, \cite{BCP}, \cite{BBBP}, \cite{BBP} on the Kac Ising
models.  
Recent developments of the theory have been prompted by the work
\cite{Cerf}  and comprise a rigorous derivation of the (free)
multi-dimensional  Wulff crystal shapes in the context of Bernoulli bond
percolation  \cite{Cerf} and in the context of nearest neighbor Ising
model in \cite{Bo} and  \cite{CePi}. Most recently the theory has been 
generalized to 
the case of symmetric $q$-states Potts models \cite{CePi2}. 
We refer to the review article  \cite{BIV}
for a detailed account of the  $\bbL_1$-theory  and its relation to the
Dobrushin-Koteck\'{y}-Shlosman approach (see e.g. \cite{DKS}, \cite{Pfister},
\cite{I1}, \cite{I2},  \cite{PV2}, \cite{IS}) to the problems of
phase segregation.

\subsection{Other microscopic derivations of the Winterbottom construction}

The derivation of the macroscopic geometry of an equilibrium droplet in contact with a
substrate has been accomplished in several models, which can be divided into two
classes.

The first class is that of SOS-type models, with a positivity constraint,
a fixed volume under field, and an attraction to the wall. The first such
study was done for the $1+1$-dimensional model in~\cite{CDR} (see also \cite{MR}
for an alternative approach), while its higher dimensional counterpart was
treated in~\cite{BI}; see also~\cite{CMR} for an analysis of the effect of
substrate roughness in 1+1 dimensions. Other related works, which neglect the
interaction between the substrate and the interface (therefore describing a completely
neutral substrate), are~\cite{BD}, \cite{DGI} and~\cite{DM}; notice
however that these works cover only one special point of the full regime
of complete wetting. This class of models provide a somehow appealing description of
this phenomenon, but has several shortcomings. The first one is that the full
phenomenology of the wetting transition cannot be observed, due to the
very particular constraint that the droplet has to be the graph of a function
(e.g. it is of course impossible for the droplet to detach from the wall). The
corresponding variational problem also displays some non-physical
properties, e.g. the minimizer are in general not scale invariant. And of course,
these SOS models are effective models for interfaces, not truly microscopic ones, so
they
do not clarify the mechanism which generates the interface itself. On the
other hand, due to these simplifying features, it is possible to obtain
some additional information about fluctuations around the equilibrium
shape, see~\cite{DGI}, which go  beyond what we can achieve for the models
considered in the present work.
For an in-depth review and references on these effective models, we
refer the reader to the survey~\cite{Gi}.

The second class contains more realistic lattice gas models of the type
we consider here. Only the two-dimensional nearest-neighbor
Ising model had been studied before~\cite{PV1,PV2}, but the description
provided is rather complete. The wetting transition can be fully
described, and the characterization of typical configurations is very accurate,
especially when these results are combined with the local limit estimates
of~\cite{IS}:
There is a unique macroscopic droplet whose boundary is close in Hausdorff
distance to one of the solutions of the Winterbottom variational problem,
and all the other droplets are of size at most $\log L$, $L$ being the linear
size of the box. Notice that at low temperature, it is also possible to obtain
precise informations on the fluctuations around the minimizer, 
see~\cite{DH}.

\subsection{Microscopic model: Bulk phases}
\label{subsection_bulk}

We start by introducing the measures describing finite-range ferromagnetic
Ising models with pair interactions in the bulk; the complete description,
including the influence of the substrate, is given in
Subsection~\ref{subsection_wallfree energy}.

The system is contained in the finite box
$$
\bbB_N \df \setof{i\in\bbZ^d}{\abs{i_k}\leq N,\, k=1,\dots, d}\,.
$$
The interactions are given by a set of real numbers ${\bf J} \df
(J_{ij})_{i,j\in\bbZ^d}$, with the following properties:
\begin{enumerate}
\item $J_{ij}\df J_{\abs{i-j}}$;
\item the graph $\left(\bbZ^d ,\cE_{\bf J}^d\right)$ with  edge set $ \cE_{\bf
J}^d$ consisting of all (unoriented) pairs of vertices $(i,j)$  with $J_{|i-j|}
>0$ is connected;
\item $J_k\geq 0$ for all $k\in\bbZ^d$;
\item $J_k = 0$ if $\normI{k}>R$, for some $R<\infty$ (the smallest such $R$ is
called the range of the interaction).
\end{enumerate}

Let $\overline\gs \in \{-1,1\}^{\bbZ^d}$. The Hamiltonian in $\bbB_N$ with
boundary condition $\overline\gs$ is given by
$$
{\bf H}^{\overline\gs}(\gs) \df 
- \frac{1}{2} \sum_{\pairs ij \subset \bbB_N} J_{ij}
\gs_i\gs_j - \sum_{i\in \bbB_N,\, j \in\bbZ^d\setminus\bbB_N} J_{ij}
\gs_i{\overline\gs}_j\,,
$$
and the corresponding (finite-volume) Gibbs measure at inverse temperature $\gb$
is the probability measure on $\{-1,1\}^{\bbZ^d}$ defined by
$$
\Is^\gb_{N,\overline\gs} (\gs) \df
        \begin{cases}
                \frac 1{{\bf Z}^{\overline\gs}_{\bbB_N}}\, \exp\left(  -\gb{\bf
                H}^{\overline\gs}(\gs)  \right)
                & \text{if $\gs_i={\overline\gs}_i$, for all
                $i\not\in\bbB_N$,}\\
                0 & \text{otherwise.}
        \end{cases}
$$
Two boundary conditions of particular importance are the $+$ b.c. defined by
$\overline\gs\equiv 1$, and the $-$ b.c. $\overline\gs\equiv-1$; the
corresponding finite-volume Gibbs states are denoted by $\Is^\gb_{N,+}$ and
$\Is^\gb_{N,-}$ respectively. It follows from correlation inequalities that the
weak limits $\mu^\gb_+ \df \lim_{N\rightarrow\infty} \mu^\gb_{N,+}$ and $\mu^\gb_-
\df \lim_{N\rightarrow\infty} \mu^\gb_{N,-}$ exist and are translation
invariant (in fact ergodic). It is a basic result of Statistical Mechanics that
in dimensions $d\geq 2$ there exists $0<\gb_{\rm c}(d,\bf J)<\infty$ such that
$\mu^\gb_+ = \mu^\gb_-$ if $\gb < \gb_{\rm c}$
(this also implies uniqueness of the infinite-volume Gibbs
state), while $\mu^\gb_+ \neq \mu^\gb_-$ if $\gb > \gb_{\rm c}$. The order
parameter associated to this phase transition is the
spontaneous magnetization,
$$
m^*(\gb) \df \lim_{N\rightarrow\infty} \mu^\gb_+\left(\frac 1{\abs{\bbB_N}}
\sum_{i\in\bbB_N} \gs_i \right) = \mu^\gb_+(\gs_0)\,.
$$
It can indeed be shown that $m^*(\gb)=0$ when $\gb < \gb_{\rm c}$ and
$m^*(\gb)>0$ when $\gb>\gb_{\rm c}$. Of course, by symmetry,
$\mu^\gb_-(\gs_0)=-m^*(\gb)$.

\subsection{Microscopic model: Surface tension}  
\label{subsection_tension}

Let $\vec n\in\bbS^{d-1}$ be a unit normal. Assume for the definiteness that 
$(\vec n,\vec{e}_d)>1/\sqrt{d}$. Given two positive real numbers $L$ and $M$,
define
\begin{equation}
\label{1_st_bc}
\Lambda_{M,L}(\vec n)~\df~\left\{\ x\in\bbR^d:\ 
\begin{split}
&\abs{(x,\vec{e}_k)} \leq L/2\ \ \text{for}\ k=1,\dots ,
d-1\\
&\text{and}\ \abs{(\vec n, x)}\leq M/2 
\end{split}
\right\}
\end{equation}
Thus,  $\Lambda_{M,L}(\vec n)$  is a parallelepiped with the base orthogonal to
$\vec n$ and having area $L^{d-1}/(\vec n,\vec{e}_d )$ and height
$M$. With a slight abuse of notation we identify  $\Lambda_{M,L}(\vec
n)$ with its intersection with $\bbZ^d;\ \Lambda_{M,L}(\vec n)
=\Lambda_{M,L}(\vec n)\cap\bbZ^d$. Let $\PF_{\Lambda_{M,L}(\vec n)}^+$ and
$\PF_{\Lambda_{M,L}(\vec n)}^{\pm}$ be the partition  functions on
$\Lambda_{M,L}(\vec n)$ with respectively  ``$+$'' and $\vec n$ boundary
conditions, the latter being defined by $\gs_i = \text{sign}((\vec n, i))$,
with $\text{sign}(0)=1$. The (per unit area) surface tension  of the 
$\pm$-interface stretched in the direction orthogonal to $\vec n$ is defined via: 
\begin{equation}
\label{1_tau_beta}
\st (\vec n)~\df~-\limtwo{L\to\infty}{M\to\infty}
\frac{(\vec n,\vec{e}_d)}{L^{d-1}}
\log\frac{\PF_{\Lambda_{M,L}(\vec n)}^{\pm}}{\PF_{\Lambda_{M,L}(\vec n)}^+}
\end{equation}
The important ferromagnetic feature of the model (or, possibly, rather of the
proof of the corresponding fact)  is that the limit in \eqref{1_tau_beta} is
defined and {\em does not} depend on the order in which the numbers 
$L$ and $M$ go to infinity. In ~\cite{miracle} this result has 
been formulated  
 for the sequences of domains tending to $\infty$ in the sense of Landford,
 but the corresponding proof goes through also in the 
$\lim_{M\to\infty} \lim_{L\to\infty}$ case. 

Finally,  as 
it has been  proven in~\cite{LebPfister}, 
 the surface tension is uniformly strictly
positive whenever $\beta>\beta_{\rm c}$.

\subsection{Microscopic model: Wall free energy}  
\label{subsection_wallfree energy}  
Microscopic models for the substrate surface have been studied in 
\cite{FroehlichPfister87a} and  \cite{FroehlichPfister87b} in the context of 
the nearest neighbor model. Let us  reformulate,  with brief comments on the
validity in the general case we consider here, all those of the results of
\cite{FroehlichPfister87b} which we are going to use in the sequel.

Define the lattice half-space $\halfspace \df \setof{i\in\bbZ^d}{i_d > 0}$. Our
spin configurations $\gs$ are the elements of $\{ -1,1\}^{\halfspace}$.  The
microscopic influence of the substrate is modeled by magnetic fields (chemical
potentials in the lattice gas language) $\bdf =(\bdf_1 ,\dots ,\bdf_r )$ acting
on spins in the first $r$ microscopic layers\footnote{We do not assume any
connection between  the number $r$ and the range of the interaction $\{
J_{ij}\}$. } of $\halfspace$. Thus, the formal  Hamiltonian on $\{
-1,1\}^{\halfspace}$ is given by:
\begin{equation}
\label{1_boundary_ham}
\Ham^{\bdf} (\sigma ) ~\df~-\sum_{ \pairs ij\subset\halfspace}J_{ij}\gs_i\gs_j ~-~
\sum_{k=1}^r\  \sum_{i\in\halfspace\, :\, i_d =k}\bdf_k\sigma_i .
\end{equation}
Given two different vectors $\bdf$ and $\bdf^{\prime}$ of boundary magnetic
fields we say that $\bdf$ is larger than $\bdf^{\prime}$ if
$\bdf_k\geq\bdf^{\prime}_k$ for all $k=1,\dots,r$.

In order to describe the  appropriate surface energy corrections induced by the
substrate, we put it in  contact with either of the two ``$-$'' or ``$+$'' bulk
phases.  This is done in a  standard way: Given $N\in\bbN$ we construct the
finite lattice box
\begin{equation*}
\BoxxN \df \setof{i\in\halfspace}{-N < i_k \leq N,\, k=1,\dots,d-1,\quad
0< i_d \leq N}\,,
\end{equation*}
and consider the models with formal Hamiltonian \eqref{1_boundary_ham} and,
respectively, with ``$-$'' and ``$+$'' boundary conditions on
$\halfspace\setminus \BoxxN$. Let $\Isbd{-}{\BoxxN}$ and $\Isbd{+}{\BoxxN}$ be
the corresponding finite volume Gibbs states on $\{-1 ,1\}^{\BoxxN}$ and 
$\PFbd{-}{\BoxxN}$, $\PFbd{+}{\BoxxN}$ be the associated partition functions.
The difference $\Delta_{\beta ,\bdf}$ between the interfacial free energies is
defined via:
\begin{equation}
\label{1_Delta_beta_eta}
\Delta_{\beta ,\bdf}~\df~\lim_{N\to\infty}\frac1{(2N)^{d-1}}\log
\frac{\PFbd{+}{\BoxxN }}{\PFbd{-}{\BoxxN}}.
\end{equation}
\begin{pro}[\cite{FroehlichPfister87b}] The limit in \eqref{1_Delta_beta_eta} 
is well defined and  monotone non-increasing  in $\bdf$.
Furthermore, uniformly in $\bdf$
\begin{equation}
\label{1_st_bound}
\left|\Delta_{\beta ,\bdf}\right|~\leq ~\st (\vec{e}_d ).
\end{equation}
\end{pro}
As in \cite{FroehlichPfister87b} the assertion of the theorem  follows from the
FKG properties of the ferromagnetic measures. We refer to Section~3.1 in
\cite{FroehlichPfister87b} (with obvious modifications to fit the general case
we consider here) for  details.
\begin{remark}
On the heuristic level the bound in \eqref{1_st_bound} should be
clear: one possible scenario under the $\Isbd{-}{\BoxxN}$ measure is to create
a microscopically divergent film of ``$+$'' phase along the  substrate surface
reducing, thus, the ratio in \eqref{1_Delta_beta_eta} to the  one appearing in
the definition of $\st (\vec{e}_d )$ in \eqref{1_tau_beta}.  This scenario
certainly becomes dominant at large positive values of $\bdf$, which happens to
be equivalent to the uniqueness of the limiting boundary  Gibbs state
\cite{FroehlichPfister87b}.  A non-trivial analysis of the corresponding phase
diagram (in terms of $\bdf$) is, however,  an almost entirely open question,
c.f. \cite{Chalker}, \cite{BDZ}, \cite{CV} for the related results for
effective interface models.
\end{remark}

\subsection{The result}
\label{subsection_result}
In principle, the results of the paper should hold for any $\beta >\beta_c$.
However, our approach relies on the validity of Pisztora's renormalization
machinery, which we shall describe in detail in
Section~\ref{section_renormalization}.
Consequently, we need the following 
 assumption on the inverse temperature $\beta$:\\[2mm]
{\bf (A)} Pisztora's coarse graining holds at the inverse temperature
$\beta$.\\[2mm]
We use $\frB$ to denote the set of $\beta>\beta_c$ for  which the above
assumption holds.\\[2mm]
Define the  total average magnetization on $\BoxxN$ as 
$$
{\bf M}_{\BoxxN}(\gs )~\df~
\frac{1}{\left| \BoxxN \right|} \sum_{i \in \BoxxN} \gs_i .
$$
 The averaged magnetization ${\bf M}_{\BoxxN}$
concentrates (as $N\to\infty$) under the measure $\Isbd{+}{\BoxxN}$ on the
spontaneous ``$+$''~phase magnetization $m^* (\beta )$. In other words the bulk
phase induced by $\Isbd{+}{\BoxxN}$ is the purely ``$+$'' one. In order to
enforce phase segregation one should introduce a canonical type constraint
which would shift the average magnetization ${\bf M}_{\BoxxN}$ inside the
phase  coexistence region. In the lattice gas language this would amount to
fixing the average number of particles strictly between the two extremal
equilibrium densities. The canonical constraints we are working with here are of
the  integral type: Given $m\in \left(-m^* (\beta ),m^*(\beta )\right)$, let us
define the canonical conditional measure
\begin{equation}
\label{1_m_measure}
\Isbd{m,+}{\BoxxN}(\gs )~=~\frac1{\PFbd{m,+}{\BoxxN}}\exp\left[
-\gb\,\Ham_\BoxxN^\bdf(\gs)\right] \1_{\{{\bf M}_{\BoxxN}\leq m\}}(\gs
)~=~\Isbd{+}{\BoxxN}\left(\sigma ~\Big|~{\bf M}_{\BoxxN}\leq m\right) .
\end{equation}
The impact of such a conditioning on the bulk properties should  amount to the
creation of a ``$-$''~phase island of macroscopic size close to 
\begin{equation}
\label{1_v_m}
v(m) ~\df ~ \frac{m^*-m}{2m^*} \, .
\end{equation}
Our microscopic justification of the Winterbottom construction 
above the complete wetting threshold; $\Delta_{\beta ,\bdf} >-\tau_\gb (\vec{e}_d )$, 
 gives a rigorous
meaning to the heuristic picture above and, moreover, asserts that under the
natural scaling the shape of the ``$-$''~phase island converges to the 
Winterbottom shape $\cK_{\gb ,\bdf}\lb v (m)\rb$. Such a result is, of course, 
impossible unless $\cK_{\gb ,\bdf}\lb v (m)\rb$ fits inside the unit box (the 
macroscopic vessel of the system)
$$
\uBox~\df~[-\frac{1}{2},\frac{1}{2}]^{d-1}\times [0,1],
$$
which is obtained by the scaling of $\BoxxN$ by the factor $1/(2N)$. Define
$$
\bar{m}(\gb ,\bdf )~=~\min\setof{\, m}{\cK_{\gb ,\bdf}\lb v (m)\rb \subseteq \uBox\,}.
$$
\begin{thm}
\label{theorem_A}
Assume that $\beta\in\frB$ and all the components of the boundary 
magnetic field 
$\bdf\in\bbR^r$ have the same
sign. Assume, furthermore, that $\Delta_{\beta ,\bdf}> - \tau_\beta (\vec{e}_d)$. 
Then,
  for every $m\in ] \bar m (\gb ,\bdf ) ,m^*(\gb) [$,
\begin{equation}
\label{claim_A}
\begin{split}
\lim_{N\to\infty}\frac1{N^{d-1}}\log\frac{\PFbd{m,+}{\BoxxN}}{\PFbd{+}{\BoxxN}}~&=~
\lim_{N\to\infty}\frac1{N^{d-1}}\log\Isbd{+}{\BoxxN}
\left( {\bf M}_{\BoxxN} <m \right)\\
&=~- \cW_{\beta ,\bdf}\left( \cK_{\beta ,\bdf }(v(m)) \right)\,
\df\, -w^*_{\gb,\bdf}(m).
\end{split}
\end{equation}
\end{thm}

Theorem~\ref{theorem_A} gives a microscopic justification of the
thermodynamical variational problem leading to the Winterbottom 
construction, 
 but, as it is, 
 yields  little information on the statistical properties of the
microscopic spin configuration under the canonical measure
$\Isbd{m,+}{\BoxxN}$. In the $\bbL_1$-approach we pursue here the microscopic
spin fields are identified only on the (renormalized) mesoscopic level through
the local order parameters or, equivalently, through locally averaged
magnetization profiles. These local order parameters are piece-wise constant
functions on the continuum box $\uBox$ and they are constructed from the
microscopic spin configuration $\gs\in\{-1,1\}^{\BoxxN}$ in the following
way:\\[2mm]
{\bf Choice of scales.} We shall always take the microscopic  size of the
binary form, $N=2^n$. Similarly,  all the intermediate mesoscopic scales are of
the  form $K=2^k ,k\in\bbN$.\\[2mm]
{\bf Partition of $\BoxxN$. } At each fixed mesoscopic scale $K=2^k$   we split
the microscopic vessel $\BoxxN$  into the disjoint union of shifts of the
mesoscopic box $\bbB_K \df \{-\tfrac12 K,\dots ,\tfrac12 K\}^d$.  These shifted
boxes are centered at the lattice points from the rescaled set  $\bbD_{N,K} \df
K\left(\bbD_{N/K}-(1/2,\dots ,1/2)\right)$:
\begin{equation}
\label{Ksplit}
\BoxxN ~=~\bigvee_{i\in\bbD_{N,K}}\bbB_K (i) ,
\end{equation}
where $\bbB_K (i)\df i+\bbB_K$.\\[2mm]
{\bf Induced partition of $\uBox$.} We scale \eqref{Ksplit} by the factor 
$N$.
With each mesoscopic lattice box $\dBox{K} (i)$ in the partition 
\eqref{Ksplit} we associate the continuum box 
$$
\sBox{N,K}(x)~ \df~ x +\left(-\tfrac12 \frac KN,\tfrac12 \frac KN\right]^d
$$ 
centered at the point $x=\frac{1}{N} i$.  We use $\sDor{N,K}$ to
denote the set of all such centers $x$; $\sDor{N,K} = \frac{1}{N} \bbD_{N/K}$. 
The induced mesoscopic partition of $\uBox$  is given by:
\begin{equation}
\label{ksplit}
\uBox~=~\text{closure}\left(\bigvee_{x\in\sDor{N,K}}\sBox{N,K} (x)\,\right) \, .
\end{equation}
\noindent
The local magnetization profile $\cM_{N,K}$ corresponding to a microscopic spin
configuration $\gs \in \{-1,1\}^{\BoxxN}$ is a function  on the continuum box
$\uBox$, piece-wise constant with respect to the partition  \eqref{ksplit}. For
every $x\in\sDor{N,K}$ the value of $\cM_{N,K}$ on $\sBox{N,K} (x)$  equals to
the averaged magnetization on the corresponding lattice box  $\dBox{K} (i)$
with $i=N x$;
$$
\cM_{N,K}(\gs)(y)~=~\frac1{K^d}\sum_{j\in\dBox{K} (i)}\gs_j \, , 
\qquad\text{for}\ \ y\in\sBox{K}(x).
$$
Informally, $\cM_{N,K}(\gs )$ is the resolution with which one observes the
microscopic spin field $\gs$. Local proximity of the Gibbs measure to one of
the pure phases is, in this way, quantified by the proximity of $\cM_{N,K}(\gs
)$ to the corresponding order parameter $\pm m^* (\beta )$. Our probabilistic
counterpart of the thermodynamic  Theorem~\ref{theorem_A} states that on
coarser scales local  magnetization profiles $\cM_{N,K}$ comply, with the
overwhelming conditional  probability $\Isbd{m,+}{\BoxxN}$, with the
thermodynamic prediction. This agreement has to be understood in the
$\bbL_1$-sense: Given a measurable subset  $A\subset\bbR^d$ let us define the
function $\1_A\in \bbL_1 (\uBox )$ via:
\begin{equation*}
\1_A (x)~=~\left\{
\begin{split}
& +1,\qquad \text{if}\ x\in A\\
& -1,\qquad\text{if}\ x\not\in A
\end{split}
\right.
\end{equation*}
\begin{thm}
\label{theorem_B}
Fix a number $\nu <1/d$ and 
assume that $\beta\in\frB$ and all the components of the boundary 
magnetic field 
$\bdf\in\bbR^r$ have the same
sign. Assume, furthermore, that $\Delta_{\beta ,\bdf}> - \tau_\beta (\vec{e}_d)$. 
 Let $m\in \, ] \bar m (\gb ,\bdf ) ,m^*(\gb) [\,$.  
 Then for every $\gd >0$ fixed, one can choose a
finite mesoscopic scale $K_0 =K_0 (\gb ,\bdf ,\gd )$, such that
\begin{eqnarray*}
\lim_{N \to \infty} \;  \min_{K_0 \leq K \leq N^\nu} \; 
\Isbd{m,+}{\BoxxN} \left(
\min_{x } 
\|\cM_{N,K}(\gs) - m^*_\gb 
\1_{\{x + \cK_{\beta ,\bdf}(v(m)) \} } 
\|_{\bbL_1 (\uBox )} \leq\gd
\right)~
 = ~ 1 \, .
\end{eqnarray*}
\end{thm}

Of course, the above minimum can be restricted to the shifts of 
$ \cK_{\beta ,\bdf}(v(m))$ along the wall $\cP_d$ in the partial wetting/drying 
case; $\left| \Delta_{\gb ,\bdf }\right| <\tau_\gb (\vec{e}_d )$, and, respectively, 
to all admissible shifts of the Wulff shape $\cK_\gb (v )$ within $\uBox$ 
in the case of complete drying; $\Delta_{\gb ,\bdf } = \tau_\gb (\vec{e}_d )$.

\subsection{Structure of the paper and further remarks}
\label{subsection_structure}
To a large extent the proof of Theorems~\ref{theorem_A} and~\ref{theorem_B} is
a book-keeping exercise based on the renormalization formalism we have
developed for the nearest-neighbor Ising model  in \cite{BIV}. In that paper
we tried to decouple deep model oriented facts from the relatively soft
$\bbL_1$ techniques. The fact that the  proof swiftly goes through in the
general case we consider here is a dividend of such an approach. 

The deep model oriented fact is the validity of Pisztora's coarse
graining~\cite{Pisztora} and its application to the relaxation properties of
the FK measures~\cite{CePi}; both have been originally developed 
in the nearest
neighbor context, but go through literally without changes in any model with
translation invariant finite-range ferromagnetic pair interactions. In
Section~\ref{section_FK} we briefly recall the  corresponding construction and 
show how to adjust it to treat the boundary  surface tension. 
In Section~\ref{section_renormalization} we set up the machinery of random 
mesoscopic phase labels and relate it to the context of the theory of 
functions of bounded variation.  Since the 
boundary field introduces a singularity of the surface tension at
the bottom wall, some 
additional care (as compared to
\cite{BIV}) is needed at this stage. Furthermore, the appropriate FK-representation 
  in the case of negative boundary magnetic
fields happens to be a conditional one.  In particular, the 
compactness estimates have to be modified.

The proof of Theorem~\ref{theorem_B} (and, on 
the way, Theorem~\ref{theorem_A})
is relegated to  Section~\ref{section_proof}. In view of the results and
adjustments of the  
preceding sections it closely follows the pattern laid down
in \cite{BIV}. Accordingly, our exposition in this section  
will be rather
concise with repeated references to the appropriate parts of
\cite{BIV}.

Finally, we would like to point out that though the minimizers of the variational problem 
$\text{{\bf (WP)}}_{\gb ,\bdf}$ restricted to the sets $P\subseteq\uBox$ are, in general,
not known for $m <\bar{m}_{\gb ,\bdf}$,  
the ${\Bbb L}_1$-approach still could be 
pushed through  
 to yield meaningful results in the sense that one can prove that the mesoscopic configurations
concentrate around a (in general unknown) set of surface energy minimizers.
 Such an idea of decoupling between the probabilistic analysis of the microscopic model and 
the investigation of the limiting variational problem has been put forward in   \cite{PV1} 
 and, recently, systematically  exploited in \cite{CePi2}. 
 Indeed, 
the microscopic derivation of phase coexistence is based upon local decoupling estimates 
 and on surgical 
procedures in small regions localized along the interface of the
crystals.
Therefore, from the probabilist point of view the  approach is local and does not rely on
the global shape of the solutions to the corresponding variational problems.  On 
the technical level, however, one should   
adjust the 
approximation results of Subsection~\ref{subsection Approximation result I} 
to the case of
unknown minimizers. The approximation procedure we use here 
hinges on the fact that the minimizers 
are known to be convex. 
If the surface tension is singular at the boundary (that is in the 
partial wetting/drying cases),  
 a complete geometric analysis without such an information 
 would require a  serious 
detour into the geometric
measure theory which, from our point of view,  would not bring any deeper insights on the 
microscopic phenomena involved in the Winterbottom construction. We refer to 
\cite{Cerf} or to \cite{CePi2} for 
a detailed exposition of the appropriate facts from the theory 
of functions of bounded variation.

A more challenging problem is to gain a better understanding of the 
statistical properties of phase boundaries, even on mesoscopic or 
macroscopic scales. For example, the only information the 
${\Bbb L}_1 $-approach yields on the 
entropically repulsed interface is  a localization on 
the macroscopic scale,  
which
is, of course, not such  a fantastic result, 
since we absolutely do not control
the typical   microscopic size of the  corresponding  fluctuations.

\section{FK Representation}
\label{section_FK}
\setcounter{equation}{0}
FK representation is an artificial coupling between  the ferromagnetic spin
models and the, so called, random cluster measures. Its power rests with the
fact that the latter are in the {\em uniqueness} regime and possess uniform 
decoupling and relaxation properties, even when the original spin system is in
the phase coexistence region\footnote{Strictly speaking, this is known to be
true for all sub-critical temperatures, except possibly for a countable subset
of the latter, see~Subsection~\ref{ssec_MPL}.}. The art of the FK
representation is, therefore, to single out such  spin configuration properties
which admit a reformulation in the FK language.

\subsection{Representation of the bulk states}
\label{subsection_bulk_fk}
With a given finite range translation invariant pair interaction potential
${\bf J}=\{J_{ij}=J_{|i- j|}\}$ we associate the  graph $\left(\bbZ^d ,\cE_{\bf
J}^d\right)$ where the  edge set $ \cE_{\bf J}^d$ consists of  all (unoriented)
pairs of vertices $(i,j)$  with $J_{|i-j|} >0$. Each such pair will be called a
bond of $ \cE_{\bf J}^d$. The set $\Omega_{{\bf J}}\df\{ 0,1\}^{\cE_{\bf
J}^d}$  is the sample space for the dependent percolation measures associated
with ${\bf J}$. Given $\go\in\Omega_{\bf J}$ and a bond $b=\pairs ij\in\cE_{\bf
J}^d$, we say  that $b$ is open if $\go (b)=1$. Two sites of $\bbZ^d$ are said
to be connected if one can be reached from another via a chain of open bonds.
Thus, each $\go\in\Omega$ splits  $\bbZ^d$ into the disjoint union of maximal
connected components, which are called the open clusters of $\Omega$. Given a
finite subset $B\subset\bbZ^d$  we use $c_B (\go )$ to denote the number of
different open finite 
clusters of $\go$ which have a non-empty intersection
with $B$.

We define next the FK measures which correspond to the finite volume spin Gibbs
states on the boxes ${\bbB}_N$: First of all these measures put non-trivial
weights on the percolation configurations  $\omega\in \Omega_{N,{\bf J}}\df \{0
,1\}^{\cE_{N,{\bf J}}^d}$,  where the set ${\cE_{N,{\bf J}}^d}$ comprises those
of the bonds $b\in\cE_{\bf J}$ which intersect ${\bbB}_N$. The boundary
conditions are specified by a frozen percolation configuration  $\pi\in
\Omega_{N,{\bf J}}^c\df \Omega_{\bf J}\setminus\Omega_{N,{\bf J}} $. Using the
shortcut  $c^\pi_N (\go ) =c_{{\bbB}_N} (\go\vee \pi )$ for the joint
configuration $\go \vee\pi\in\cE_{\bf J}$, we define the finite volume FK
measure $\FKm{\beta ,\pi}{N}$ on $\Omega_{N,{\bf J}}$  with the boundary
conditions $\pi$ as:
\begin{equation}
\label{FKm}
\FKm{\beta ,\pi}{N}\left(\go \right)~\df~\frac1{\FKpf{\beta ,\pi}{N}}
\left( \prod_{b\in
\cE_{N,{\bf J}}} \big( 1-p_{\gb}(b) \big)^{1-\go_b}
\big( p_{\gb}(b) \big)^{\go_b} 
\right)\, 2^{c^{\pi}_N(\go)}\,,
\end{equation}
where, for a bond $b=\pairs ij\in \cE_{N,{\bf J}}$, we define the  corresponding
percolation probability $p_{\gb}(b) =1-\exp(-2\gb J_{|i-j|})$.

The measures $\FKm{\beta ,\pi}{N}$ are FKG partially ordered with respect to
the lexicographical order of the boundary condition $\pi$. Thus, the extremal
ones correspond to the free ($\pi\equiv 0$) and wired ($\pi\equiv 1$) boundary
conditions and are denoted as  $\FKm{\beta ,{\rm f}}{N}$ and $\FKm{\beta ,{\rm
w}}{N}$  respectively. The  corresponding infinite volume ($N\to\infty$)
limits  $\FKm{\beta ,{\rm f}}{\,}$ and $\FKm{\beta ,{\rm w}}{\,}$ always exist,
and  it takes a relatively soft argument \cite{Grimmett} to show that for all
but at most a countable set of temperatures $\beta\in\bbR_+$ there is a unique
infinite volume FK-measure:\\[2mm]
{\bf (A1)} \hspace*{1cm} We assume that the inverse temperature $\beta$
satisfies $\FKm{\beta ,{\rm f}}{\,}=\FKm{\beta ,{\rm w}}{\,}$. \\[2mm]
Assumption {\bf (A)} is a combination of {\bf (A1)} and the condition 
that there
is percolation in slabs, see Assumption {\bf (A2)} below.  

The spontaneous magnetization $m^* (\beta )$ admits the following expression in
FK terms:
\begin{equation}
\label{FK_m_star}
m^* (\beta )~=~\lim_{N\to\infty}\FKm{\beta ,{\rm w}}{N}\left( 0\lra
{\bbB}_N^c\right)~=~ \FKm{\beta ,{\rm w}}{\,}\left( 0\lra\infty \right).
\end{equation}
More generally,  the finite volume (spin) Gibbs state $\IsN$ on $\{
-1,1\}^{{\bbB}_N}$ (see Subsection~\ref{subsection_bulk}) can be recovered
\cite{EdwardsSokal}  from the wired FK-state $\FKm{\beta ,{\rm w}}{N}$ as
follows: Sample a  percolation configuration $\omega\in \Omega_{N,{\bf J}}$ from
$\FKm{\beta ,{\rm w}}{N}$. Spins at the sites of ${\bbB}_N$ which belong to the
wired component of  $\omega\vee{\bf 1}$ (that is connected to the boundary
${\bbB}_N^c$) are assigned  value $+1$, whereas the remaining open clusters of
$\omega\vee{\bf 1}$ are painted to $\pm 1$ with probability $1/2$ each.

\subsection{Representation of the surface tension}
\label{subsection_tension_fk}
The surface tension of the spin model \eqref{1_tau_beta} can be expressed in
the FK language in the following way:  Let us split the  boundary $\partial
\Lambda_{M,L}(\vec n) = \bbZ^d\setminus  \Lambda_{M,L}(\vec n)$ into two
pieces:
\begin{align*}
\partial \Lambda_{M,L}(\vec n)~&=~\partial^+ \Lambda_{M,L}(\vec n)
\bigcup \partial^- \Lambda_{M,L}(\vec n)\\&\df~
\left(\partial \Lambda_{M,L}(\vec n)\cap
\setof{i\in\bbZ^d}{i_d >0}\right)\bigcup\left(
\partial \Lambda_{M,L}(\vec n)\cap\setof{i\in\bbZ^d}{i_d \leq 0}\right)
\end{align*}
Then the ratio of partition functions in \eqref{1_tau_beta} is equal to
\begin{equation}
\label{1_tau_beta_fk}
\frac{\PF_{\Lambda_{M,L}(\vec n)}^{\pm}}{\PF_{\Lambda_{M,L}(\vec n)}^+}~=~
\FKm{\beta ,{\rm w}}{\Lambda_{M,L}(\vec n)}
\left(\partial^+ \Lambda_{M,L}(\vec n)
\nlra
\partial^- \Lambda_{M,L}(\vec n)
\right),
\end{equation}
where the FK-percolation event  $\left\{\partial^+ \Lambda_{M,L}(\vec n) \nlra
\partial^- \Lambda_{M,L}(\vec n) \right\}$ means that no site of
$\Lambda_{M,L}(\vec n)$ is connected both  to $\partial^+ \Lambda_{M,L}(\vec
n)$ and to $\partial^- \Lambda_{M,L}(\vec n)$.\\

The expression~\eqref{1_tau_beta_fk} for the surface tension is not very
convenient in practice. Indeed,
in higher dimensions one should not expect the interface to decouple
along  
pure
wired b.c. It would therefore be very useful to have a more robust
definition which would still fit in the framework of the 
$\bbL_1$-theory.   The idea is twofold. First, one can relax the
pinning of the ``microscopic interface'' along the lateral sides of
$\Lambda_{M,L}(\vec n)$, by requesting only that the top and bottom faces of
the box are disconnected. Of course in order to still recover the surface
tension in the direction orthogonal to $\vec n$, one has to choose $M\ll L$.
The control of the interactions on the lateral sides is one of the
major technical contribution initiated by Cerf in \cite{Cerf}.
The second point is to relax boundary conditions, i.e. to replace the wired
b.c. by some arbitrary $\pi$. In order to avoid a possible deformation of the
``microscopic interface'' due to the attraction by some b.c. $\pi$, one has to
let enough room for the system to relax to its unique equilibrium phase. This
is done by imposing the boundary condition $\pi$ outside the bigger box
$\Lambda_{L,L}$.

More precisely, let $\partial^{\rm top} \Lambda_{M,L}(\vec n)$, resp.
$\partial^{\rm bot} \Lambda_{M,L}(\vec n)$, be the face of $\partial^+
\Lambda_{M,L}(\vec n)$, resp. $\partial^- \Lambda_{M,L}(\vec n)$, normal to
$\vec n$. As a consequence of the relaxation properties of FK measures
derived by Cerf and Pisztora~\cite{CePi} (Proposition 3.1), the following holds
\begin{lem}
\label{lem surface tension} 
Let $\gep>0$. For any $\gb\in\frB$, 
$$
\tau_\gb(\vec n) = - \frac {(\vec n, \vec{e}_d)}{L^{d-1}} \log
\FKm{\beta ,\pi}{\Lambda_{L,L}(\vec n)}
\left(\partial^{\rm top} \Lambda_{\gep L,L}(\vec n)
\nlra
\partial^{\rm bot} \Lambda_{\gep L,L}(\vec n)
\right) + c_{\gep,L}(\pi,\vec n)\,,
$$
with $c_{\gep,L}(\pi,\vec n)$ going to zero uniformly in $\pi$ and $\vec n$, as
$L$ goes to infinity and $\gep$ goes to zero.
\end{lem}

\subsection{Representation of the boundary states} 
The FK-notation for the boundary Gibbs states closely follows the full space
setup introduced in the previous subsection. The relevant graph for the
inter-spin  interactions is $\left(\halfspace ,\cL_{\bf J}^d\right)$ which is
just the  restriction of $\left(\bbZ^d ,\cE_{\bf J}^d\right)$ to the half-space
$\halfspace$. In order to incorporate the boundary magnetic field $\bdf
=(\bdf_1 ,...,\bdf_r )$ (see \eqref{1_boundary_ham}) we augment this graph with
a ghost site $\frg$  connected to all the sites in the first $r$ layers of
$\halfspace$. Thus the  edge set for the boundary model is given by
$$
\cL_{{\bf J},\bdf}^d~\df ~\cL_{{\bf J}}^d\, \bigcup\,\left\{ (i,\frg )~\big|~
i\in\halfspace\ \text{and}\ i_d\leq r\right\} .
$$
Similarly, given $N\in\bbN$, we use $\cL_{N,{\bf J},\bdf}^d $ to denote the 
set of bonds of $\cL_{{\bf J},\bdf}^d$ which  have a non-empty intersection 
with  $\BoxxN$. The sample space for finite volume FK states on $\BoxxN$ is given
by 
$$
\Xi_{N,{\bf J},\bdf}\df\{ 0,1\}^{\cL_{N,{\bf J},\bdf}^d}.
$$
Assume that $\bdf\geq 0$ (that is $\bdf_k\geq 0$ for $k=1,\dots ,r$). Given  a
configuration $\pi\in\{0,1\}^{ \cL_{{\bf J},\bdf}^d \setminus \cL_{N,{\bf
J},\bdf}^d}$, the  FK measure $\FKm{\gb ,\bdf}{N}$ on $\Xi_{N,{\bf
J},\bdf}$ is defined by
\begin{equation}
\label{FKm_boundary}
\FKm{\gb ,\bdf,\pi}{N}\left( \xi\right)~=~\frac1{\PF^{\gb ,\bdf}_{N}}
\left( \prod_{b\in
\cL_{N,{\bf J},\bdf}^d} \big( 1-p_{\gb}(b) \big)^{1-\xi_b}
\big( p_{\gb}(b) \big)^{\xi_b} 
\right)\, 2^{c_N^\pi (\xi)}\,,
\end{equation}
where $c_N^\pi (\xi )$ denotes the number of open finite 
clusters of $\xi\vee\pi$ 
which intersect $\BoxxN$ and do not contain the ghost site $\frg$,
 whereas the percolation probabilities $p_\gb (b)$ are
defined exactly as in  \eqref{FKm} for the bonds $b\in\cL^d$ and equal to
$1-{\rm e}^{-2\gb\bdf_k}$ on  the ghost bonds $b =(i,\frg )$ with $i_d =k$,
$k=1,\dots ,r$. We suppress the super-index $\pi$ for the wired state with the
boundary conditions $\pi\equiv 1$.

The boundary Gibbs state $\Isbd{+}{\BoxxN}$ can be reconstructed as follows: 
sample a bond configuration $\xi\in\Xi_{N,{\bf J},\bdf}$ from  $\FKm{\gb
,\bdf}{N}$, and paint with $1$ all the clusters of $\xi$ connected 
either to $\halfspace\setminus\BoxxN$ or to $\frg$, whereas all the  remaining
clusters of $\xi$ are to be painted into $\pm 1$ with probability  $1/2$ each. 
The corresponding 
joint bond-spin probability measure is denoted by $\Joint_{N}^{\gb,\bdf}$.
\begin{remark}
The FK measures $\FKm{\gb ,\bdf ,\pi}{N}$ corresponding
to  the non-negative magnetic fields $\bdf\geq 0$ are the basic ones and we
shall use them in all our  considerations.
\end{remark}

The  representation derived in the case of non-positive magnetic fields  $\bdf
= (-\abs{\bdf_1},\dots ,-\abs{\bdf_r})$ can be described as follows: Define the  FK
percolation event $\frJ_N\subset  \Xi_{N,{\bf J},\bdf}$ as
\begin{equation}
\label{J_event}
\frJ_N~\df ~\left\{ \xi\in\Xi_{N,{\bf J},\bdf}~\Big|~\frg\nlra \halfspace\setminus
\BoxxN\right\},
\end{equation}
and set $\FKm{\gb ,\bdf}{N}\left(\,\cdot\,\right)=
\FKm{\gb ,\abs{\bdf} }{N}\left(\,\cdot\,\big|\,\frJ_N\right)$. Then, 
the boundary Gibbs measure $\mu_{\BoxxN ,+}^{\gb,\bdf}$ can be reconstructed from
$\FKm{\gb ,\bdf}{N}\left(\,\cdot\,\right)$ as above, except that the 
spins connected to $\frg$ are, this time, painted into $-1$. 
The corresponding 
joint measure is denoted by $\jpm$. With a slight abuse of notation
we shall write (keeping in mind that the $\frg$-cluster is 
repainted into $-1$):
\begin{equation}
\label{JN_conditioning}
\jpm\left(~\cdot ~\right)~=~\Joint_{N}^{\gb,|\bdf |}
\left(~\cdot~\Big|\frJ_N\right) .
\end{equation}
%
\subsection{Representation of the wall free energy}
\label{subsection_wallfreeenergy_fk}
In view of the spin-flip symmetries the difference $\Delta_{\gb ,\bdf}$ (see 
\eqref{1_Delta_beta_eta}) admits the following expression in FK terms
(recall the convention $\bdf\geq 0$):
\begin{equation}
\label{1_Delta_beta_eta_fk}
\Delta_{\gb,\bdf}~=~-\lim_{N\to\infty}\frac1{(2N)^{d-1}}
\log\FKm{\gb ,\bdf }{N}\left(
\frJ_N\right),
\end{equation}
and, accordingly, if $\eta$ is a negative boundary field
$\Delta_{\gb ,\bdf}=-\Delta_{\gb ,|\bdf|}$.\\

As in the case of the surface tension, the above definition is too restrictive
to be used in practice, because of the pure wired b.c.. Fortunately, a more
robust expression is also available here.

Let $\partial^{\rm top,\gep}\, \BoxxN \df \BoxxN \cap\setof{i\in\bbL^d}{i_d =
[\gep N]}$, where $[x]$ denotes the integer part of $x\in\bbR^d$. Then, using
the relaxation properties of FK measures (\cite{CePi}, Proposition~3.1),
one can prove the following
\begin{lem}
\label{lem wall free energy}
Let $\gep>0$. For any $\gb\in\frB$ and any $\bdf\geq 0$,
$$
\Delta_{\gb,\bdf}~=~-\frac1{(2N)^{d-1}}
\log\FKm{\gb ,\bdf, \pi}{N}\left( ~\frg\nlra \partial^{\rm top,\gep}\, \BoxxN
\right) + c_{\gep,N}(\pi)\,,
$$
with $c_{\gep,N}(\pi)$ going to zero uniformly in $\pi$, as N goes to infinity
and $\gep$ goes to zero.
\end{lem}

\section{Renormalization}
\label{section_renormalization}
\setcounter{equation}{0}

In order to perform the analysis of phase coexistence in a macroscopic
setting, a convenient formalism, the Geometric Measure Theory, is introduced 
in subsection \ref{ssec_BV}.
The embedding of the discrete spin system in the continuum setting relies
on the introduction of renormalized variables, the mesoscopic phase labels
(Subsection \ref{ssec_MPL}).

\subsection{Functions of bounded variation}
\label{ssec_BV}

We refer the reader to \cite{EG} for an introduction to the theory 
of functions of bounded variation and
to \cite{BIV} for a discussion on its relevance in the context of
phase coexistence.

Let $\cO$ be an open smooth neighborhood of $\uBox$.
On the macroscopic scale the system is characterized by
the boundary condition $g\in{\rm BV}\left( \cO\setminus\uBox ,\{\pm 1\}
\right)$ 
and by  a $\pm 1$ phase function $u\in\BV$. 
The fact that $u (x) = 1$ for some $x\in{\rm int}\uBox$ means that
locally at  $x$ the system is in equilibrium in the ``$+$'' phase.
Thus, $u$ should be interpreted as a signed indicator of
the regions containing the different phases and the boundary of 
the set $\{ u\vee g = -1\}$, where
\begin{equation}
\label{patched_ug}
u\vee g (x)~=~
\left\{
\begin{split}
\ u(x)\qquad &{\rm if}\ x\in{\rm int}\uBox\\
\ g(x)\qquad &{\rm if}\ x\in \cO\setminus\uBox
\end{split}
\right.
\end{equation}  
as the interface in the presence of $g$-boundary condition.

It is well known \cite{EG} that $u\vee g\in{\rm BV}(\cO ,\{\pm 1\})$ whenever 
 the phase function $u\in\BV$. 
Now, for any $v$ in ${\rm BV}(\cO,\{\pm 1\})$, 
there exists a generalized notion of the boundary of 
$\{ v = - 1 \}$
called reduced boundary and denoted by $\partial^* v$.
If $\{ v = - 1 \}$ is a regular set, $\partial^* v$ coincides with the usual 
boundary $\partial v$. Given a boundary condition $g\in {\rm BV}(\cO\setminus
\uBox , \{\pm 1\})$ and a phase function $u\in\BV$ we use $\partial_g^* u$ to 
denote the reduced boundary of $u$ in the presence of the b.c. $g$:
\begin{equation}
\label{partial_ug}
\partial_g^* u~=~\partial^*(u\vee g)\cap\uBox~=~
\partial^*(u\vee g)\setminus \partial^* g .
\end{equation}
In the setting we are working with here, there 
are two natural types of boundary conditions 
depending on the sign of the magnetic field either
$g=g_+\equiv 1$ in the case of positive magnetic fields $\eta\geq 0$, or
\begin{equation}
\label{g_pm}
g(x)~=~g_\pm (x)~=~\left\{
\begin{split}
\ +1,\qquad &{\rm if}\ x_d > 0\\
\ -1,\qquad &{\rm if}\ x_d \leq 0 
\end{split}
\right.
\end{equation}
in the case of negative magnetic fields $\eta < 0$.

For any non-negative boundary magnetic field $\eta$  and any 
$g\in{\rm BV}\left( \cO\setminus\uBox ,\{\pm 1\}
\right)$ define the functional 
$$
\intSTbdf{\bdf} \lb u\,\big|\, g\rb ~=~
\int_{\partial^*_g u \setminus \cP_d} \tau_\gb (\vec{n}_x) d \cH^{(d-1)}_x
+ \taubd  \int_{\partial^*_g u \cap \cP_d}  \; d \cH^{(d-1)}_x,
 .
$$
Set $\cW^{+}_{\gb,\bdf }(\,\cdot\,) = \intSTbdf{\bdf} \lb \cdot\,\big|\, g_+\rb$ and, 
accordingly, $\cW^{\pm}_{\gb,\bdf }(\,\cdot\,) = \intSTbdf{\bdf} \lb \cdot\,\big|\, g_\pm\rb$.
 As we shall see below this notation exactly corresponds to the 
color blind formalism of the FK-representation and to the construction of
 the measures $\jpm$ in \eqref{JN_conditioning}. Thus, 
the Winterbottom functional \eqref{intro_functional}
can be rewritten (recall that for $\bdf \leq 0$, 
$\taubd =-\Delta_{\gb, \abs\bdf}$)  as
\begin{equation}
\label{functional_g}
\intSTbdf{\bdf} (u) ~=~\left\{
\begin{split}
 &\cW^{+}_{\gb,\bdf }(\, u\,)\ \ \qquad\ \ \ \ \ \ \   \text{for}\ \bdf\geq 0\\ 
&\cW^{\pm}_{\gb,|\bdf | }(\, u\,) -\Delta_{\gb,\abs \bdf}
\ \ \ \ \text{for}\ \bdf\leq 0
\end{split}
\right.
\end{equation}
%
%
%
In this way the functional $\intSTbdf{\bdf}$ is unambiguously defined on
$\BV$ for every constant sign magnetic field $\bdf$.  
%
 
For $m\in ]\hat{m} (\gb ,\bdf ),m^* (\beta ) [$ the Winterbottom shape  
$\cK_{\gb ,\eta} (v(m))$ fits into
 $\uBox$ and the function $\1_{\cK_{\gb ,\eta} (v(m))}\in \BV$ is  the
stable minimizer of $\cW_{\gb ,\bdf}$ in the following sense: For every $\delta>0$
there exists $\epsilon >0$ such that for every function 
 $u\in\BV$ with 
$\int_{\uBox}u(x)\,{\rm d}x\leq m/m^*_\beta$, 
\begin{equation}
\label{stability}
\cW_{\gb ,\bdf} (u)\,\leq \,w^*_{\gb ,\bdf}(m) +\epsilon\ \Longrightarrow\ 
\min_x \Big\| u\,-\,\1_{x+\cK_{\gb ,\eta} (v(m))}\Big\|_{{\Bbb L}_1 (\uBox )}\,
\leq \, \delta .
\end{equation}
Indeed, let $\tau_{\gb ,\bdf}$ be the support function of the Winterbottom 
shape $\cK_{\gb ,\bdf}$ itself. Regardless of the sign of the monochromatic 
magnetic field $\bdf$ define the relaxed functional: 
$$
\widehat{\cW}_{\gb ,\bdf} (u) ~=~\int_{\partial_{g_+}^* u}
\tau_{\gb ,\bdf}(\vec{n}_x)d\cH_x^{d-1} .
$$
  Being the Wulff type 
functional (or, in an alternative 
 terminology, being the mixed volume \cite{Sch})
 $\widehat{\cW}_{\gb ,\bdf}$ is 
lower-semicontinuous with respect to the ${\Bbb L}_1 (\uBox )$-convergence.
Furthermore, by the refinement \cite{FM} of 
the Brunn-Minkowski inequality \cite{Sch}, the shifts
and dilations of $\cK_{\gb ,\bdf}$ are the only minimizers of 
$\widehat{\cW}_{\gb ,\bdf}$ in the problems  with the corresponding volume 
constraints. Hence, $\widehat{\cW}_{\gb ,\bdf}$ satisfies the stability
property \eqref{stability}.
 However, $\cW_{\gb ,\bdf} \geq \widehat{\cW}_{\gb ,\bdf}$ and their values 
coincide on the shifts and dilations of the Winterbottom shape 
$\cK_{\gb ,\bdf}$. As a result,  the stability property \eqref{stability} 
is immediately inherited by the original functional $\cW_{\gb ,\bdf}$.

 In the partially wetting case $\Delta_{\beta ,\eta}  < \tau_\beta
( \vec{e}_d )$ the minimum in \eqref{stability} is over all 
admissible shifts $x\in \cP_d$ along the substrate wall. In order to see
 this notice that in the latter 
situation the minimum $w^*_{\beta ,\eta }(m )$ is
strictly less than the unconstrained minimum of the Wulff functional 
$\cW_\beta (u)$ over the functions $u\in \BV$;
 $\int_{\uBox}u(x)\,{\rm d}x\leq m/m^*_\beta$. 

Fix now  a point 
$x\in\uBox\setminus\cP_d$ and define $r=x_d >0$. For any $\delta >0$ and 
any function $u\in\BV$ satisfying 
$\| u-\1_{x+\cK_{\gb ,\eta} (v(m))}\|_{{\Bbb L}_1 (\uBox )}\,
\leq \delta$   one choose a section $r_\delta \leq r$ such that 
$\cH^{d-1} \lb\{x_d=r_\delta\}\cap\setof{x}{u(x)=-1}\rb\leq \delta/r$. 
Modifying $u$  as $u_\delta =+1$ on $x_d <r_\delta$ we obtain that 
$\cW_\beta (u_\delta )\leq \cW_{\beta ,\eta}(u ) +\delta \max\tau_\beta/r$ 
 and, of course, that 
 $\int_{\uBox}u_\delta (x)\,{\rm d}x\leq (m+\delta )/m^*_\beta$. As a result,
 a $\cW_{\beta ,\eta}$-minimizing sequence cannot converge to 
 $\1_{x+\cK_{\gb ,\eta}}$.

\subsection{Mesoscopic phase labels}
\label{ssec_MPL}

In the $\bbL_1$-approach the 
local order parameter $\cM_{N,K}$ is quantified by three  
values $\{-1,0,1\}$ according to the
local proximity of the system to one of the pure phase. Such a 
renormalization
procedure is delicate to implement directly. One way, proposed by Pisztora in
\cite{Pisztora},  relies on the FK representation of Ising model. A technical
advantage to work with the FK representation, is the uniqueness of 
the measure in
the thermodynamic limit even when there is breaking of symmetry  for the Ising
model. In a sense the FK measure is much less sensitive to boundary effects
than the Gibbs measure, which underlines the required decoupling properties.\\

\noindent
\underline{\it Typical configurations for the FK representation :}

As explained in \eqref{FK_m_star}, the FK counterpart of spontaneous
magnetization is a uniform positive probability that the site 0 is connected to
the boundary of arbitrarily large boxes with wired boundary conditions. The
implementation of Pisztora's coarse graining requires an {\it a priori} stronger
notion, namely percolation in slabs :
%
\\[2mm]
{\bf (A2)} \hspace*{1cm}
$\exists \gd > 0, \exists L_0 > 0, \forall L > L_0, \qquad
\lim_{N \to \infty} \ \inf_{x,y \in S_{L,N}}
\FKm{\gb, \rm f}{S_{L,N}} ( x \lra y) > \gd \, ,$
\\[2mm]
where $S_{L,N}$ is the slab  $\{ \abs{i_1} \leq L, -N \leq i_k \leq N, \ k = 2,
\dots , d \}$. 
The importance of this notion was first realized in \cite{ACCFR} in
the context of Bernoulli percolation. 
If {\bf (A2)} is satisfied then there is percolation
and $\gb > \gb_c$. The critical value $\hat \gb_c$  above which {\bf (A2)}
holds is called slab percolation threshold.
We define $\frB$ as the subset of $]\hat \gb_c,\infty[$ for which assumption
{\bf (A1)} holds. The set $\frB$ differs from $]\hat \gb_c,\infty[$ by at most
a countable number of points. In fact, it is conjectured that $\frB =
]\gb_c,\infty[$.

Originally (see \cite{Pisztora}) the estimates on typical 
configurations were devised in a nearest neighbor set-up. 
Inspection of the proof shows that the results extend readily to
models with finite range interaction.

\medskip

Property 
{\bf (A2)} enables to obtain enhanced estimates on the percolation
in domains by slicing these domains into slabs. Furthermore assumption {\bf
(A1)} implies a relaxation property of the system  and therefore uniform
estimates can be deduced for arbitrary boundary conditions.

We can now proceed to describe Pisztora's coarse graining. Let $\ga$ be in 
$(0,1)$ and fix $\gz>0$. On a coarse grained scale $K \in \bbN$, the typical
configurations $\go$ in the box $\dBox{K}$ satisfy the 3 properties below :
\begin{enumerate}
\item There is a unique crossing cluster $C^*$ in  $\dBox{K}$, i.e. a cluster
which is connected to all the faces of the vertex boundary  
$$
\partial \dBox{K}
\df \{ x \in \dBox{K}^c \ | \ \exists y \in  \dBox{K}, \text{ such that } 
J(x-y) >0
\}.
$$
The faces are simply the faces of the box $\dBox{K}$ enlarged by a factor
depending on the range of the interaction. 
\item Every open path in $\dBox{K}$  with diameter  larger than $K^\ga$ is
contained in $C^*$. (The diameter of a subset $A$ of  $\bbZ^d$ is $\sup_{x,y
\in A} \|x - y\|_1$.)
\item The density of the crossing cluster in $\dBox{K/2}$ is close to $m^*$  with
accuracy $\gz >0$
\begin{eqnarray*}
| C^* | \in [m^* -\gz, m^* + \gz] \, \frac{K^d}{2^d},
\end{eqnarray*}
where $| \cdot |$ denotes the number of vertices in a set.
\end{enumerate}

\smallskip

A configuration $\go$ in $\dBox{K}$ is good if it satisfies  the 3 assertions
above. It was proven in \cite{Pisztora} that for any $\gb \in \cB$ and  for
large enough mesoscopic scales, good configurations are, uniformly over 
boundary conditions, typical in the following sense :
%
\begin{eqnarray}
\label{typical}
\inf_{\pi} \  \FKm{\gb, \pi}{K} 
\left( \go \text{ is a good configuration in} \  \dBox{K} \right)
\geq 1 - \exp ( - c_\gz K^\ga) \, ,
\end{eqnarray}
where $c_\gz > 0$ depends on the accuracy $\gz$ of the coarse graining.

\medskip

\noindent
\underline{\it Mesoscopic phase labels :}


On the mesoscopic scale $K = 2^k$, each box $\sBox{N,K}(x)$, centered at $x \in
\sDor{N,K}$ is labeled by the variable $u^{\gz ,{\bf FK}}_{N,K} (\go,x)$
\begin{equation*}
\forall x \in \sDor{N,K}, \qquad
 u^{\gz,{\bf FK}}_{N,K} (\go,x)  
\df
\begin{cases}
1 &     \text{if } \go \text{ is a good configuration in } 
\dBox{2 K} (Nx),\\
0 &     \text{otherwise,}
\end{cases}
\end{equation*}
where the information recorded in $u^{\gz,{\bf FK}}_{N,K} (\go,x)$ takes into
account the bond configuration in the larger box $\dBox{2 K} (Nx)$. We say that
a block $\sBox{N,K}(x)$ is regular if $u^{\gz,{\bf FK}}_{N,K}(x) =1$. 

Estimate \eqref{typical} is uniform over the boundary conditions, therefore it
enables to control the variable $u^{\gz,{\bf FK}}_{N,K}(x)$  
independently of the
events which occurs on boxes which are not  $*$-connected to $\sBox{N,K}(x)$.
Following \cite{Pisztora}, \cite{LSS}, this implies stochastic domination of
the the field $\{ x\in
\sDor{N,K}~|~u^{\gz,{\bf FK}}_{N,K} (\go,x)=0 \}$ by the 
Bernoulli site percolation  measure with parameter  $\exp ( -
c_\gz K^\ga)$ on $\sDor{N,K}$; in particular,
\begin{equation}
\label{domination1}
\RCbdf{\wired}{\abs\bdf} \left(  u^{\gz,{\bf FK}}_{N,K} (x_1 )=0, \dots ,
u^{\gz,{\bf FK}}_{N,K} (x_\ell )=0\right)~\leq \exp ( - c_\gz K^\ga \ell) \, .
\end{equation}

Furthermore, the basic FK-estimate \eqref{typical} enables to  
control the local 
magnetization profiles $\cM_{N,K}$,
for large enough mesoscopic scale $K=2^k$,  by a coarse graining which
satisfies similar decoupling properties.  Let $x$ be in $\sDor{N,K}$. Then for
any $y \in \sBox{N,K} (x)$, we define the mesoscopic phase label $u^\gz_{N,K}$
as: 
\begin{equation}
\label{eq_magnetization}
u^\gz_{N,K} (\gs,\go,y) \df
\begin{cases}
{\rm{sign}}(C^*)  &     \text{if } u^{\gz,{\bf FK}}_{N,K} (\go,x)  = 1
\text{ and } \abs{\cM_{N,K} (\gs, x) - {\rm{sign}}(C^*) \,  m^*} 
< 2 \gz,\\
0 & \text{otherwise,}
\end{cases}
\end{equation}
where $C^*$ is the crossing cluster in $\dBox{K} ( x)$.\\

A block $\sBox{N,K} (x)$ may have label 0 for two reasons : either 
$u^{\gz,{\bf FK}}_{N,K} (\go,x) = 0$, or $\go$ 
is a good configuration on $\bbB_{2K} (Nx)$, but  the
local average magnetization $\cM_{N,K}(x)$ on $\bbB_{K}(Nx)$ is non typical.
The former case is covered by \eqref{domination1}. 
In the latter case, the corresponding shift of the magnetization is 
entirely due to an abnormal coloring of the small (i.e. with the diameter
less than $K^\alpha$ ) $\omega$-clusters intersecting $\bbB_{K}(Nx)$. 
 The probability of
such a deviation is, for each particular block $\bbB_K (Nx)$, bounded above
by ${\rm exp}\{ - c\gz^2 K^{d(1-\alpha )}\}$. 
On the other hand, all small
clusters intersecting $\bbB_K (Nx)$ 
lie inside $\bbB_{2K}(Nx)$, which already leads to the 
required decoupling properties. 
Notice that the coloring argument inside a good block is 
insensitive to the repainting of the $\frg$-cluster into $-1$ on the 
event $\frI_N$. 
We refer to the original article 
\cite{Pisztora} for the precise workout of these estimates and to
\cite{Bo} where the proof  has been adapted to the 
setup we employ here.

As a result we obtain that there there exists a sequence $\{\gr^\gz_K\}$ with
$\lim_{K\to\infty}\gr^\gz_K =0$, such that the field
$\{ x \in \sDor{N,K}~|~ u^\gz_{N,K} (\gs,\go,x ) =0\}$ is stochastically
dominated by the Bernoulli site percolation  measure  $\bbP_{\text{\rm
perc}}^{\gr^\gz_K}$ on $\sDor{N,K}$. In particular, for $\bdf \geq 0$
\begin{equation}
\label{domination2}
\jp \left( u^\gz_{N,K} (x_1 )=0,...,u^\gz_{N,K} (x_\ell )=0\right)~\leq 
(\gr_K^\gz)^\ell \, .
\end{equation}
%
%
\begin{remark}
\label{0-block}
If for two different points $x,y\in \sDor{N,K}$ the corresponding
$u^\gz_{N,K}$-phase labels have opposite signs, that is if 
$u^\gz_{N,K} (x) u^\gz_{N,K} (y)=-1$, then
on any finer scale $K^{\prime} \leq K$ any $*$-connected chain of 
$\sBox{N,K^{\prime}}$ blocks joining $\sBox{N,K}(x)$ to $\sBox{N,K}(y)$ 
contains at least one block with zero $u^{\gz,{\bf FK}}_{N,K^{\prime}}$
(and, hence, with zero $u^{\gz}_{N,K^{\prime}}$)-label.
This follows from the fact that 2 boxes with opposite labels cannot have
a common crossing FK-cluster;
thus they are necessarily separated by a contour of 0-blocks.
\end{remark}

\begin{remark}
The modification of the boundary interaction has no impact on the coarse
graining.  
Since the slab percolation threshold hypothesis requires percolation on
slab with free boundary condition, the presence of a positive  
boundary magnetic 
field can only improve the estimates in view of FKG monotonicity.
Furthermore, in each box touching the wall, the event that a crossing cluster
is connected to the ghost point $\frg$ occurs with overwhelming probability 
(the slab percolation hypothesis {\bf (A2)} is valid uniformly over the 
boundary conditions). 
Finally the control of the density of the crossing cluster (property 3 in the
definition of a good block) is 
unchanged because the boundary field has no influence on the bulk properties
due to the hypothesis ${\bf (A1)}$.
\end{remark}

\medskip

The mesoscopic phase labels $u_{N,K}^\gz$ are $\{\pm 1,0\}$-valued function on 
$\uBox$ which are piece-wise constant with respect to the partition 
\eqref{Ksplit}.
%
%
%
%
We are now going to prove that the $\bbL_1$-difference between the local 
magnetization profiles $\cM_{N,K}$ and the  phase labels $u^\gz_{N,K}$ 
vanishes on the exponential
scale.  Thus, from the point of view of the $\bbL_1$-theory these objects
are indistinguishable, and one is entitled to switch freely from one to 
the other. 
 In particular, the fundamental exponential tightness result which we shall
 establish below in Subsection~\ref{Tightness theorem}
 comes naturally in the context of the $\{\pm 1, 0\}$-valued phase labels.

%

\subsection{Relation to magnetization profiles}
%
\begin{lem}
\label{lem coarse graining}
Let a positive 
boundary magnetic field $\eta \geq 0$ and a number $\nu <1/d$ be fixed.
For any $\gd > 0$, one can choose the accuracy $\gz$ of the coarse graining 
and a 
finite scale $K_0 =K_0 (\gd,\gb )$,
such that 
%
\begin{equation}
\label{Mkukbound}
\frac{1}{N^{d-1}}
\max_{K_0\leq K \leq N^\nu} \log \bbP_N^{\gb ,\eta } \left( 
\| \cM_{N,K}  -m^{*} u^\gz_{N,K} \|_{\bbL_1} >\gd
\right)~\leq ~ -c (\gz, K)  \; N^{1 -d \nu} .
\end{equation}
\end{lem}
\begin{remark}
By \eqref{1_Delta_beta_eta_fk} and \eqref{JN_conditioning} the 
super surface order exponential
estimate \eqref{Mkukbound}
holds, in the case of negative boundary magnetic field $\eta\leq 0$, also for 
the conditional measure $\jpm$.
\end{remark}
 
\begin{proof}
%
%
%
%
Recall that
the coarse graining is initially 
defined in terms of the mesoscopic  FK-variables 
$u^{\gz,{\bf FK}}_{N,K} (\go)$ which depend
 only on the  bond configurations. The phase label 
$u^\gz_{N,K} (\gs,\go)$ takes also into account the random coloring of
good FK-blocks. In fact, 
$$
\| \cM_{N,K}  -m^{*} u^\gz_{N,K} \|_1~\leq ~ \gz +\frac{2}{|\sDor{N,K} |}
\sum_{x\in\sDor{N,K}} 1_{u^{\gz}_{N,K} (x)=0} ,
$$
and the claim of the lemma follows by the domination by Bernoulli 
site percolation \eqref{domination2}, actually  for any choice 
of $\gz <\gd$.
\end{proof}

\subsection{Tightness theorem for mesoscopic phase labels}
\label{Tightness theorem}

We shall formulate a  tightness result which holds under fairly 
general 
phase boundary conditions outside $\uBox$. The proof is a straightforward
adaptation of the argument used to verify the claim of the tightness 
Theorem~2.2.1 in \cite{BIV} and, as we shall briefly indicate in 
Remark~\ref{rem_q} below, it can be trivially extended to cover the 
case of $q$-valued phase labels as, for example, recently considered in
\cite{CePi2}. 

Given $N=2^n$ and $K=2^k$ ($n\geq k$) consider,
as in \eqref{Ksplit}, the following splitting of $\bbR^d$:
\begin{equation}
 \label{Rd_split}
\bbR^d ~=~\bigvee_{x\in\bbZ_{N,K}^d}\sBox{N,K} (x),
\end{equation}
In what follows we shall use the same notation as in Subsection~\ref{ssec_BV}.
Thus, $\cO$ is a fixed smooth neighborhood of $\uBox$.  A boundary 
condition is an element $g\in {\rm BV}\left(\cO\setminus\uBox ,\{\pm1\}
\right)$, which is piece-wise constant with respect to the partition
 \eqref{Rd_split} on some $(M,L)$-scale. Notice that such boundary 
condition $g$ is automatically piece-wise constant on all further 
$(N,K)$ scales with $N\geq M$ and $N/K
\geq M/L$. 

A family $\{ u_{N,K}\}$ of mesoscopic phase labels is said to be compatible 
with a boundary condition $g$; $u_{N,K}\sim g$, if 
$g$ is piece-wise constant on the $(N,K)$-scale and the following 
 phase rigidity condition {\bf (R)} is satisfied:
\vskip 0.2cm
\noindent
{\bf (R)}\hspace*{1cm}
 Define $w_{N,K}= u_{N,K}\vee g$ as in \eqref{patched_ug}. 
Then for every $x\in \Boxx{N,K}$ and $y\in\bbZ^d_{N,K}$ with 
$w_{N,K}(x)w_{N,K} (y)=-1$ and for every $L\leq K$ any $*$-connected chain
of $\sBox{N,L}$-blocks  leading from $\sBox{N,K}(x)$ to $\sBox{N,K}(y)$
 necessarily contains at least one $\sBox{N,L}$-block with the 
$u_{N,L}$-label zero.
\vskip 0.2cm

In order to state the tightness result we need to define the perimeter of a
BV-function $u\in\BV$ under a boundary condition $g$: 
$$
{\cP}_g(u)~=~{\cH}_{d-1}\left(\partial_g^* u\right).
$$
It is well known \cite{EG} that for every $a\geq 0$ and for every boundary
condition $g$ the set
$$
\cC_g (a)~=~\left\{ u \in\BV~\Big|~ \cP_g (u)\leq a\right\}
$$
is compact in $\bbL_1 (\uBox )$.
Let $P_g$ be the minimal perimeter of the functions compatible with
the boundary condition $g$.\\

Before deriving an exponential tightness for the mesoscopic phase
labels under measures with boundary fields, we first state a general
result for an abstract measure $\bbP$ and general boundary conditions

\begin{thm}
\label{tightness}
Let a sequence of non-negative numbers $\{\rho_K \}$ satisfy 
$\lim_{K \to \infty} \rho_K =0$. Assume that a family of random mesoscopic
phase labels $\{u_{N,K}\}$ on $\uBox$ satisfies the following 
decoupling condition {\bf (D)}:
\vskip 0.2cm
\noindent
{\bf (D)}\hspace*{0.6cm} 
For every $N$ and $K$ the zero-label field 
$\left\{ x\in\bbZ^d_{N,K}  ~\big|~ u_{N,K} (x)=0 \right\}$ 
under the measure $\bbP$
is stochastically dominated by the Bernoulli site percolation 
measure $\bbP^{\rho_K}_{{\rm perc}}$.
\vskip 0.2cm

Fix a number $\nu <1/d$. Then for every $\delta >0$  and for each 
boundary condition $g$ there exists a 
finite scale $K_0 =K_0 (d,\delta )$ and a constant $c = c(d)>0$, such
that, uniformly in $a > P_g$,
$$
\limsup_{N\to\infty}\max_{K_0\leq K\leq N^\nu}
\frac1{N^{d-1}} \log \bbP \left[ u_{N,K}\sim g~;~
u_{N,K}\not\in \cV\left(\cC_g (a ),\gd\right)\right]~\leq~
-\frac{c}{K_0^{d-1}}a \, ,
$$
where $\cV\left(\cC_g (a ),\gd\right)$ denotes the $\gd$-neighborhood
in $\bbL_1$ of $\cC_g (a)$.
\end{thm}

\begin{remark}
\label{rem_q}
The tightness result above 
could be trivially extended to the case of $\{0 ,1,...,q\}$-valued
phase labels. One needs only to modify the rigidity condition {\bf (R)}
(to require a zero block in any $*$-connected chain between two blocks 
with different non-zero labels), and redefine $\cP_g (u)$ as the 
$\cH^{(d-1)}$-measure of the jump set of $u\vee g$~\cite{AmbrosioBraides}.
If $u_{N,K}$ is a $\{ 0,1...,q\}$-valued phase label and $g$ is a
boundary condition which satisfies
 such modified assumption {\bf (R)} and the decoupling assumption {\bf (D)},
then, for every $i=1,...,q$, the $\{ 0,\pm1\}$-label $u^i_{N,K}$ defined
as 
\begin{equation*}
u_{N,K}^i(x)~=~
\begin{cases}
1\qquad &{\rm if}\ u_{N,K}(x)=i\\
0&{\rm if}\ u_{N,K}(x)=0\\
-1&{\rm otherwise}
\end{cases}
\end{equation*}
and the boundary condition $g^i$ similarly defined 
already satisfies the assumptions of Theorem~\ref{tightness}. It remains 
to notice that
$$
\left\{ u_{N,K}\in\cV \left(\cC_g (a),\gd\right)\right\}
 \ \supset\ 
\bigcap_{i=1}^q
\left\{ u_{N,K}^i\in\cV \left(\cC_{g^i} 
\left(a/q\right),\gd/q\right)\right\}.
$$
\end{remark}

\begin{proof} The proof essentially follows Theorem~2.2.1 in 
 \cite{BIV}. The only 
difference with the periodic case discussed there is the existence
of open contours attached to the boundary. But, for 
any choice of the macroscopic boundary
condition $g$ (as defined above), these boundary contours are, in the 
language of \cite{BIV}, large, that is their diameter exceeds $C_d\log N$
on all sufficiently big $(N,K)$ renormalization scales. Let us briefly
recall the three main steps of the argument:

The first step is to check that the volume of $0$-blocks is negligible.
This simply follows by the assumption {\bf (D)} :
\begin{eqnarray}
\label{control 0 blocks}
\bbP\left( \#\{x\in\sDor{N,K}:u_{N,K} (x)=0\}\geq \gd
\left(\frac{N}{K}\right)^d
\right)
\leq \exp \left( 
- \gd  \left(\frac{N}{K}\right)^d
\log \left( \frac{\gd}{\rho_K} \right)
\right) \, .
\end{eqnarray}

The second and the third steps are
 devoted to the control of the regions surrounded by
contours, i.e. by connected surfaces of $0$-blocks.
By phase rigidity assumption {\bf (R)}
the contours play the role of mesoscopic interfaces separating regions
with mesoscopic phase labels of different signs so that they contribute
to the total perimeter of the configurations $u_{N,K}$. 
We distinguish between two types of contours, the small contours with
diameter  smaller than $C_d \log N$ (where $C_d$ is a given constant) and
the remaining contours, namely the large ones.

Peierls type estimates can be applied to bound the probability of events
such  that the collection $\{\gG_i\}$ of large contours has a total area larger
than $a (\frac{N}{K})^{d-1}$
\begin{eqnarray}
\label{control large contours}
\bbP\left( \sum_{\gG_i \ \text{large}} |\gG_i| \geq
\left(\frac{N}{K}\right)^{d-1} a \right)
\leq \exp \left( 
- c_1 a  \left(\frac{N}{K}\right)^{d-1}
\right) \, .
\end{eqnarray}

The small contours are not controlled in terms of their total area but 
of the volume of the regions they surround. By the choice of the 
macroscopic boundary condition $g$ there are no open small contours, and,
proceeding exactly as in the proof of Theorem~2.2.1 in \cite{BIV}, 
one can show that
\begin{eqnarray}
\label{control small contours}
\bbP \left( 
\sum_{\gG_i \ \text{small}} \text{vol} (\gG_i) \geq
\gd  \frac{N^d}{K^d} 
\right)
\leq \exp 
\left( - c_2   \frac{\gd N^d}{K^d \log N} \right) \, .
\end{eqnarray}
As a consequence of the previous estimate, the small contours will have
no contribution to the $\bbL_1$-norm.
Combining estimates \eqref{control 0 blocks}, \eqref{control large contours}
and \eqref{control small contours}, we 
arrive to the claim of the Theorem. 
\end{proof}

\begin{remark} The last step of the proof should be modified for 
general, that is not necessarily piece-wise constant with respect to the 
partition \eqref{Rd_split}, boundary conditions $g\in{\rm BV}\left(
\cO\setminus \uBox ,\{\pm 1\}\right)$. See also \cite{CePi} and 
\cite{CePi2} where  
small contours have been ignored and, accordingly  
tightness estimates have been derived only on   
diverging $(N,K)$-scales.
\end{remark}

Let us go back to the phase labels $\{ u_{N,K}^\gz\}$ we consider here. The 
boundary condition $g$ is given by $g=g^+\equiv 1$ in the case of 
positive boundary magnetic fields $\eta\geq 0$ and by $g=g^\pm$ specified
in \eqref{g_pm} in the case of negative $\eta < 0$.

\begin{pro}
\label{prop 1}
Let $\gb\in\frB$. Then there exists $C(\gb , d) >0$ such 
that for all $\gd$  positive,
one can find $K_0 = K_0(\gd,\gb ,d)$, $\gz = \gz (\gd)$
such that
\begin{enumerate}
\item In the case of positive magnetic fields $\eta \geq 0$,
\begin{eqnarray*}
\limsup_{N \to \infty} \;   \frac{1}{N^{d-1}}
\max_{K_0 \leq K \leq N^\nu}
\log 
\bbP_N^{\gb ,\eta} 
\left[ 
u_{N,K}^\gz  \not\in \cV\left(\cC_{g^+}(a),\gd\right)
\right] 
\leq - \frac{C}{K_0^{d-1}} a . 
\end{eqnarray*}
\item In the case of negative magnetic fields $\eta <0$,
\begin{eqnarray*}
\limsup_{N \to \infty} \;   \frac{1}{N^{d-1}}
\max_{K_0 \leq K \leq N^\nu}
\log 
\jpm 
\left[
u_{N,K}^\gz  \not\in \cV\left(\cC_{g^\pm}(a),\gd\right)
\right] 
\leq - \frac{C}{K_0^{d-1}} a + \left|\taubd\right| \, .
\end{eqnarray*}
\end{enumerate}
\end{pro}

\begin{proof}

Let us briefly work out the $\eta < 0$ case: 
Recall that $ \bbP_{N,\pm}^{\gb |\eta |} (~\cdot~) = 
\bbP_{N}^{\gb, |\eta |} (~\cdot~;~\frI_N)\big/\bbP_{N}^{\gb |\eta |}
 (~\frI_N~)$
On the other hand, any phase label $u_{N,K}^\gz$ which is compatible
 with $\frI_N$ is, after the repainting of the $\frg$-cluster into 
$-1$,  automatically compatible with the boundary condition
$g^\pm$ and, therefore, satisfies the conditions of Theorem~\ref{tightness}
(under the measure $\bbP_{N}^{\gb, |\eta |})$. Thus, the claim
 of 
the proposition
follows from the general tightness Theorem~\ref{tightness} and, in the case
of negative magnetic fields $\eta < 0$,   
 by the representation formula for the boundary surface 
tension \eqref{1_Delta_beta_eta_fk}.
\end{proof}
\vskip 0.1in

\section{The proofs}
\label{section_proof}
\setcounter{equation}{0}

We going are to derive sharp asymptotics involving the surface tension.
As mentioned in Subsection \ref{Tightness theorem}, the case of negative boundary
field requires more care, therefore we focus on this case in the following.

\subsection{Approximation result I}
\label{subsection Approximation result I}

\begin{figure}[h]
\centerline{\psfig{file=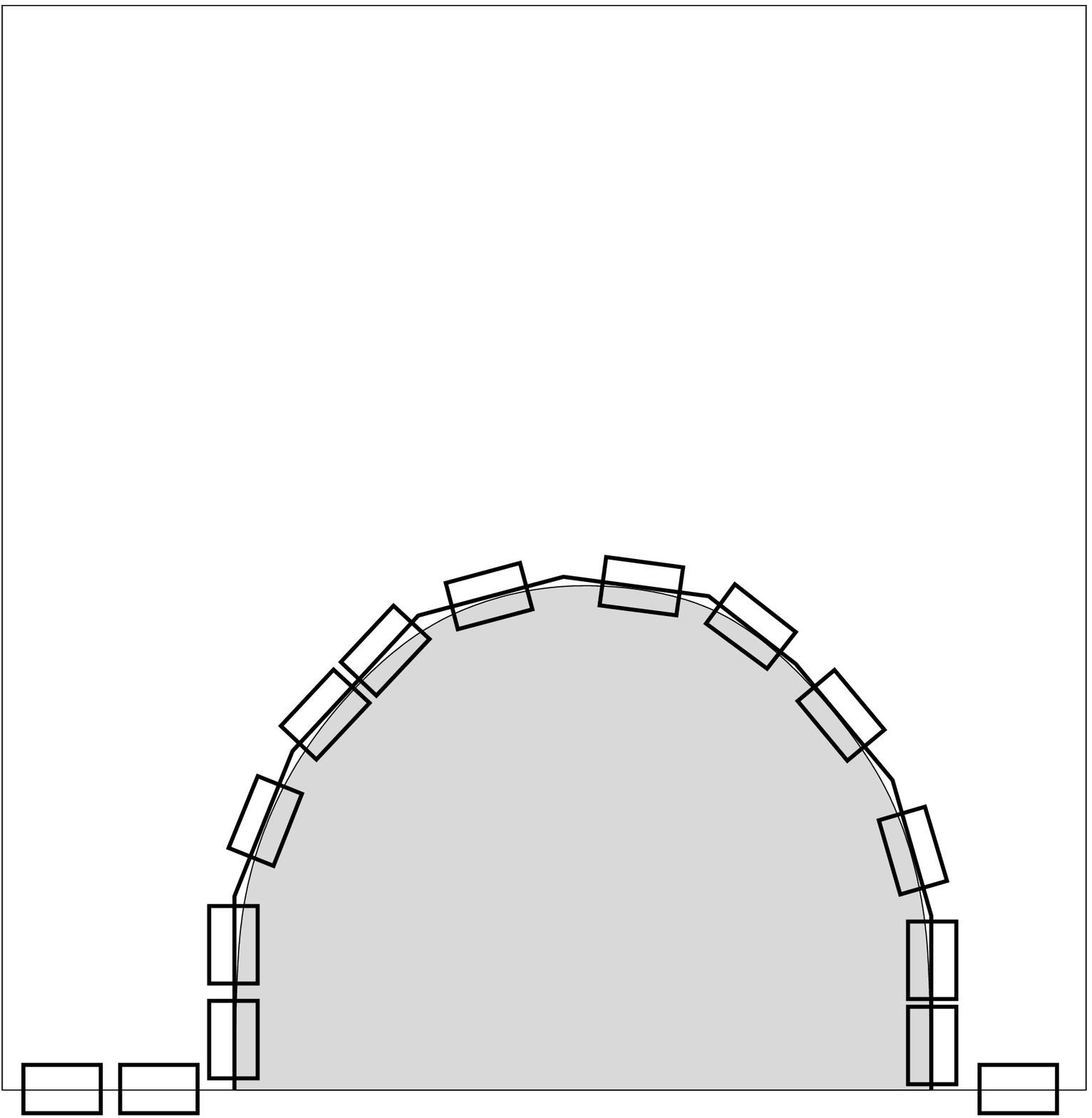,height=6cm}\hspace*{2cm}
\psfig{file=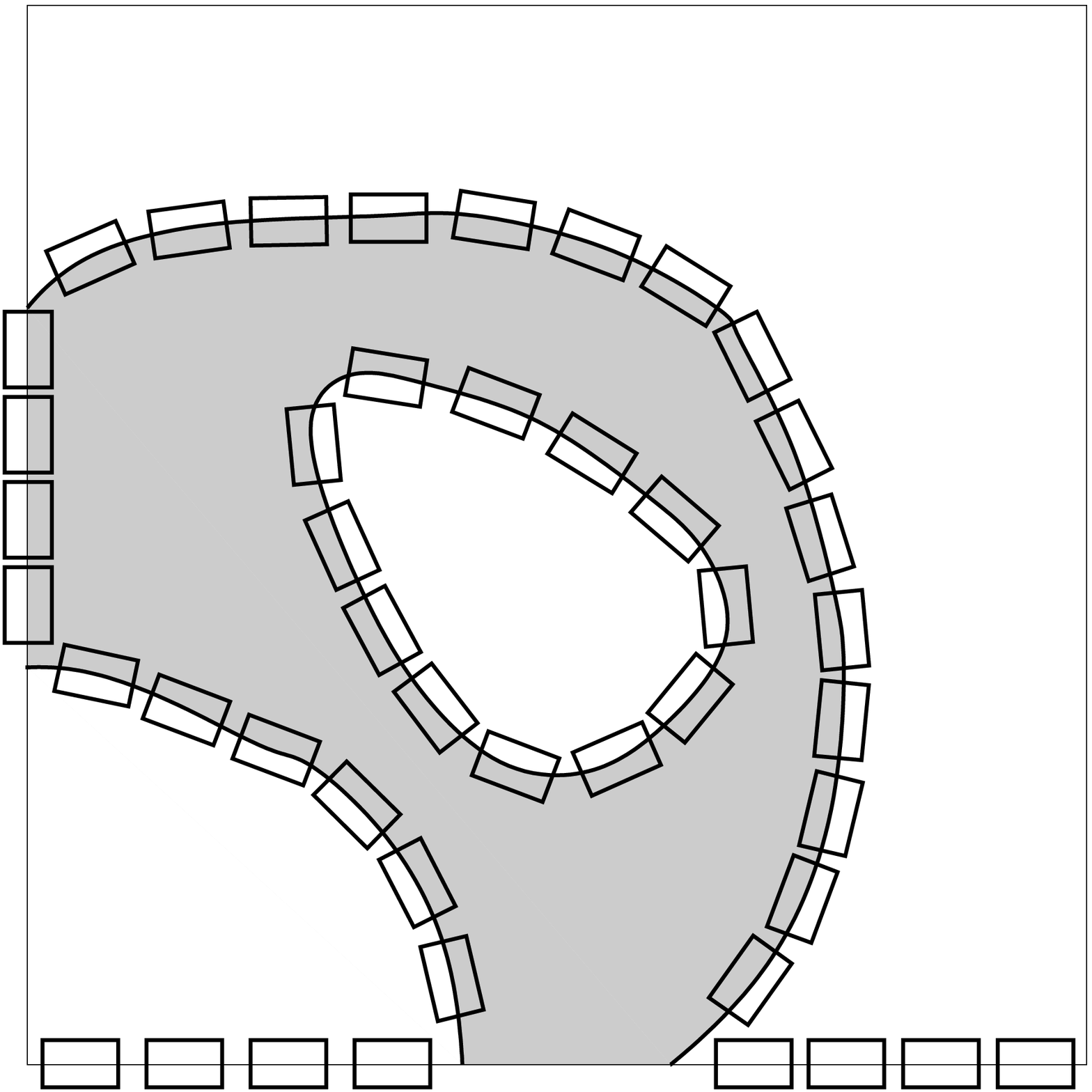,height=6cm}}
\vspace*{5mm}
\caption{Left: The approximate Winterbottom shape $\cK_{\beta ,\bdf}^\gd(v(m))$ 
and the boxes of Proposition~~\ref{prop 2}. Right:  A function $v\in \BV$ (the
shaded area corresponds to  $v\equiv -1$), and the boxes of Theorem~\ref{theo
ABCP}.}
\label{fig_approx}
\end{figure} 

Given a boundary magnetic field $\bdf$ and 
 $m \in ] \bar m_{\gb ,\bdf}, m^*_\gb[$, we are going to
approximate the Winterbottom shape $\cK_{\beta ,\bdf}(v(m))$ by regular sets
(see figure \ref{fig_approx}).

\begin{lem}
\label{thm approx}
For any $\gd>0$  one can construct  a polyhedral set 
$\cK_{\gb ,\bdf}^\gd$
 satisfying:
$$
\| \1_{\cK_{\gb ,\bdf}^\gd } - \1_{\cK_{\gb ,\bdf}} \|_{\bbL_1} \leq  \gd
\qquad {\rm  and} \qquad
\left| \intSTbdf{\bdf} \big( \cK_{\gb ,\bdf}^\gd \big) - 
\intSTbdf{\bdf} \big( \cK_{\gb ,\bdf} \big)\right|~\leq ~\gd .
$$
\end{lem}
\begin{proof}
It is a standard result of the theory of convex bodies (see e.g. \cite{Sch}~
Theorem~1.8.13) that for every $\gd >0$ there exists a convex 
polyhedral set $\cK_\gb^\gd$ such that the Hausdorff 
distance ${\rm d}_{{\rm Hausd}}( \cK_\gb^\gd ,\cK_\gb ) \leq \gd$.
Of course, the support function $\tau_\gb^\gd$
 of such $\cK_\gb^\gd$ satisfies
\begin{equation}
\label{delta_sfn}
\max_{\vec{n} \in {\Bbb S}^{d-1}}\left| \tau_\gb^\gd (\vec{n}) - 
\tau_\gb (\vec{n})\right| ~\leq ~
\gd .
\end{equation}
 Define now 
$$
\cK_{\gb ,\bdf}^\gd =\cK_{\gb}^\gd\cap
\setof{x}{x_d\geq -\taubd} .
$$ 
A comparison with \eqref{0_shape} reveals 
that 
\begin{equation}
\label{sfn_approx}
{\rm d}_{{\rm Hausd}}\big( \cK_{\gb,\bdf}^\gd ,\cK_{\gb,\bdf}\big)~\leq ~\gd\ \ 
\text{and}\ \ 
\max_{\vec{n}\in {\Bbb S}^{d-1}}\big| \tau_{\gb ,\bdf}^\gd (\vec{n}) - \tau_{\gb ,\bdf}
 (\vec{n})\big| ~\leq ~
\gd ,
\end{equation}
where  $\tau_{\gb ,\bdf}^\gd$ and $\tau_{\gb ,\bdf}$ are the support functions of 
$\cK_{\gb ,\bdf}^\gd$ and $\cK_{\gb ,\bdf}$ respectively. 
Now, for every bounded
convex $\cK\subset {\Bbb R}^d$ the volume 
${\rm vol}_d \lb \cK\rb$ can be written in 
terms of its support function $\tau_{\cK}$ as 
$$
 {\rm vol}_d \lb \cK\rb\, =\, \frac1{d}\,\int_{\partial \cK} \tau_{\cK} 
(\vec{n}_x) \, {\rm d} \cH^{(d-1)}_x . 
$$ 
Consequently, \eqref{sfn_approx} implies that
$$
\Big| \int_{\partial \cK_{\gb,\bdf}^\gd} \tau_{\gb,\bdf}^\gd
(\vec{n}_x)\, {\rm d} \cH^{(d-1)}_x \,
-\, 
\int_{\partial \cK_{\gb,\bdf}} \tau_{\gb,\bdf} (\vec{n}_x)
\, {\rm d} \cH^{(d-1)}_x \Big| ~\leq ~
d\,\gd\left| \partial \cK_{\gb,\bdf}\right| .
$$
As we have already mentioned in the end of Subsection~\ref{ssec_BV},
$$
\int_{\partial \cK_{\gb,\bdf}} \tau_{\gb,\bdf} (\vec{n}_x)
\, {\rm d} \cH^{(d-1)}_x ~=~
\cW_{\gb,\bdf}\left( \cK_{\gb,\bdf}\right) .
$$
On the other hand, for $\vec{n} = -\vec{e}_d$,  the support function
  $\tau_{\gb,\bdf}^\gd$  equals to $\taubd$.  Moreover, for every 
$x\in  \partial \cK_{\gb,\bdf}^\gd\setminus\setof{y}{y_d =-\taubd}$ any
supporting hyperplane to $\partial \cK_{\gb ,\bdf}^\gd$  at $x$ is also a 
supporting hyperplane to $\partial \cK_{\gb}^\gd$. It follows that 
 the support 
function $\tau_{\gb,\bdf}^\gd (\vec{n}_x ) =\tau_{\gb}^\gd (\vec{n}_x )$ 
 on $x\in  \partial \cK_{\gb,\bdf}^\gd\setminus\setof{y}{y_d =-\taubd}$.  
Thus, as it follows now from \eqref{delta_sfn},
$$
\left| \int_{\partial \cK_{\gb,\bdf}^\gd} \tau_{\gb,\bdf}^\gd
  (\vec{n}_x)\, {\rm d} \cH^{(d-1)}_x
-~\cW_{\gb,\bdf}\left( \cK_{\gb,\bdf}^\gd\right) \right| ~\leq ~\gd\left|
\partial\cK_{\gb,\bdf}^\gd \right| .
$$
Since by \eqref{sfn_approx}, 
 $\big|{\rm vol}_d ( \cK_{\gb,\bdf}^\gd)-{\rm vol}_d ( \cK_{\gb,\bdf})\big|\leq 
\gd\big|
\partial\cK_{\gb,\bdf}\big| $, and both 
$\big|
\partial\cK_{\gb,\bdf}\big| $ and 
$\big|
\partial\cK_{\gb,\bdf}^\gd\big| $ are bounded above  by some finite constant
$c_1 (\beta ,d )<\infty$, 
 we finally obtain the following estimate
$$
\Big| \intSTbdf{\bdf} \big( \cK_{\gb ,\bdf}^\gd \big) - 
\intSTbdf{\bdf} \big( \cK_{\gb ,\bdf} \big)\Big|~\leq ~c_2 (\beta ,d)\gd ,
$$
and the claim of the theorem follows once we redefine 
$ \cK_{\gb ,\bdf}^\gd =\cK_{\gb ,\bdf}^{\gd /(c_2\vee 1)}$.
\end{proof}
 
%
%
%
%
%
%
\subsection{The lower bound}
 As usual, we fix a number $0<\nu <1/d$.

\begin{pro}
\label{prop 2}
Assume that $\beta\in\frB$ and all the components of the boundary 
magnetic field 
$\bdf\in\bbR^r$ have the same
sign. Assume, furthermore, that 
$\Delta_{\beta ,\bdf}> - \tau_\beta (\vec{e}_d)$. 
 For every  $\gd>0$, there is a finite scale $K_0(\gd )$ such that
for every $m\in ] \bar m (\gb ,\bdf ) ,m^*(\gb) [$,
\begin{eqnarray*}
\liminf_{N \to \infty} \;   \frac{1}{N^{d-1}}
\min_{K_0(\gd) \leq K \leq N^\nu}
\log \Isbdf{+}{N}{\bdf} 
\left( \| \cM_{N,K}- m^* \1_{\cK_{\gb ,\bdf}(v\left( m\right))}
\|_{\bbL_1} \leq \gd \right) 
\geq - w_{\gb ,\bdf}^\star (m) - o(\gd) \, ,
\end{eqnarray*}
where minimal surface energy value $w_{\gb ,\bdf}^\star (m)$ has been 
defined in \eqref{claim_A}, and 
 the function $o(\cdot)$  depends only on  $\gb$ and $\eta$ and vanishes 
as the resolution $\gd$ goes to zero.
\end{pro}


\begin{proof}

We shall give the proof only in the more difficult case of negative 
boundary fields $\eta <0$.
 Also according to 
Lemma \ref{lem coarse graining}, it will be sufficient to derive
the proposition with $\cM_{N,K}$ replaced by $m^* u^\gz_{N,K}$ (for
$\gz$ small enough).

Starting from the approximate shapes $\cK_{\gb ,\eta}^\gd$ let 
 us use the transformation \eqref{1_shape} to  define the 
scaled polyhedral approximation
 $\cK_{\gb ,\eta}^\gd \left( v(m)\right) $  of 
$\cK_{\gb ,\eta}\left( v(m)\right) $.
By the approximation  result of Lemma~\ref{thm approx} it
 will  be enough to prove Proposition~\ref{prop 2} with 
$\cK_{\gb ,\eta}^\gd\left( v(m)\right) $ instead of $\cK_{\gb ,\eta}\left( v(m)\right) $. 

Define $\cP_d^{\uBox} = \cP_d\cap\uBox$. 
Since $\cK_{\gb ,\eta}^\gd\left( v(m)\right) $ is polyhedral and convex it is an 
easy matter to show that there exists a side length $h> 0$ and a finite 
number $\ell$ of disjoint parallelepipeds  
$\widehat
R^1, \dots, \widehat R^{\ell}$   with bases  
$\widehat B^1, \dots, \widehat
B^{\ell}$ included  in 
$\partial^\pm \cK_{\gb ,\bdf}^\gd \equiv \partial \cK_{\gb ,\bdf}^\gd
\symdiff \cP_d^{\uBox}$
of side length $h$ and height $\gd h$ such that:

\noindent
{\bf I\,a)}\  The sets $\widehat B^1,
\dots, \widehat B^{\ell}$ cover  $\partial^\pm \cK_{\gb ,\bdf}^\gd$ up to a set  of
measure less than $\gd$ denoted by  $\widehat U^\gd =\partial^\pm \cK_{\gb ,\bdf}^\gd
\setminus \bigcup_{i =1}^\ell \widehat B^i$  and they satisfy
\begin{eqnarray}
\label{W_pm_approx}
\Big| \sum_{\vec{n}_i\neq -\vec{e}_d}\tau_\gb (\vec{n}_i) |\widehat B^i |
\, +\,
\Delta_{\gb ,|\bdf |} \sum_{\vec{n}_i= -\vec{e}_d} |\widehat B^i |\, 
 - \cW_{\gb,| \bdf |}^{\pm}\left( \cK_{\gb ,\bdf}^\gd\right)  \Big| \leq \gd ,
\end{eqnarray}
where $\vec{n}_i$ is the unit normal to $\widehat B^i$ (and, hence, to the 
corresponding facet of $\partial^\pm \cK_{\gb ,\bdf}^\gd$). 

\noindent
{\bf I\,b)}\  For any $i=1,...,\ell$  the base
 $\widehat B^i$ divides $\widehat R^i$ into two 
parallelepipeds  $\widehat R^{i,+}$ and $\widehat R^{i,-}$, such 
that in the case of
$\vec{n}_i\neq -\vec{e}_d$,\   $\widehat R^{i,-}\subseteq \cK_{\gb ,\bdf}^\gd$ and \ 
$\widehat R^{i,+}\subseteq \uBox\setminus \cK_{\gb ,\bdf}^\gd$, whereas in the case of 
  $\vec{n}_i = -\vec{e}_d$ the corresponding base $\widehat B^i\subset \cP_d^{\uBox}$, 
 while 
$\widehat R^{i,+}\subseteq \uBox\setminus \cK_{\gb ,\bdf}^\gd$ and 
$\widehat R^{i,-}\subset \setof{x}{x_d\leq 0}$.

In order to enforce a microscopic interface close to the polyhedral set 
$K_{\gb ,\bdf}^\gd $, we define
\begin{equation}
\label{calA}
\cA = \bigcap_{i = 1}^\ell 
\{ \partial^- R^i_N \nlra \partial^+ R^i_N \}.
\end{equation}
Let us also introduce the set ${\tilde \cA}$ chosen such that for any 
configuration $\xi$ in $\cA \cap {\tilde \cA}$ the bonds inside and outside 
$N \cK_{\gb ,\bdf}^\gd$ are decoupled.
For this,
it is enough to close the bonds which are in a neighborhood of 
$U^\gd_N = N \widehat U^\gd$ (the microscopic counterpart of $\widehat U^\gd$)
\begin{eqnarray*}
\tilde \cA = \big\{ \xi \in \Xi_{N,{\bf J},\bdf}, \qquad 
\xi_b = 0, \quad  \text{if} \quad \text{dist}(b, U^\gd_N) \leq 10 \big\}\, .
\end{eqnarray*}
We get
\begin{eqnarray}
\jpm \big( \|u^\gz_{N,K}- \1_{\cK_{\gb ,\bdf}}\|_{\bbL_1} \leq \gd \big) \geq \, 
\jpm \big( \|u^\gz_{N,K}- \1_{ \cK_{\gb ,\bdf}^\gd}\|_{\bbL_1} \leq \frac{\gd}{2}  
\, | \, 
\cA\cap{\tilde \cA} \big)
\; \jpm \big( \cA\cap{\tilde \cA}\big)\, .
\end{eqnarray}

We first check that
\begin{equation}
\label{equilibre}
\liminf_{N\ra\infty} \min_{K_0(\gd)\leq K \leq N^\nu} \jpm
\big( \|u^\gz_{N,K}- \1_{ \cK_{\gb ,\bdf}^\gd}\|_{\bbL_1} \leq \frac{\gd}{2}  
\,|\,\cA\cap{\tilde \cA} \big) 
\geq \tfrac 14\, .
\end{equation}

From the argument developed in the proof of Proposition~\ref{prop 1}, we 
are going to check that 
in each of the two connected components of $\uBox \setminus \cK_{\gb ,\bdf}^\gd$
 the phase labels  $u^\gz_{N,K}$ are, 
 with a uniformly positive probability,  
 close to the equilibrium values (for the $\bbL_1$-norm).

Since $\cA \cap {\tilde \cA}$ decouples the connected components, it
is enough to consider the mesoscopic phase labels in the interior of 
the regions decoupled by $\cA \cap {\tilde \cA}$.
As the boundary of $\partial^{\pm}\cK_{\gb ,\bdf}^\gd$ is regular, 
these regions are 
fairly exhausted by the mesoscopic boxes for any scale 
$K \in [K_0,N^\nu]$.
Using the terminology of Proposition~\ref{prop 1}, we assert that
with large probability,
the mesoscopic label configurations are almost uniformly constant and contain
only contours which have no contribution in terms of the $\bbL_1$-norm.
According to estimate \eqref{control small contours}, the volume of the regions
surrounded by the small contours is arbitrarily small in the thermodynamic
limit.
Furthermore, from the usual Peierls estimates, the volume of the large contours
is negligible as well.
This implies \eqref{equilibre}.

\medskip

Finally as
$\cA\cap{\tilde \cA} \subset \frI$, we see by applying \eqref{equilibre} that
\begin{align}
\jpm(\|u^\gz_{N,K}- \1_{\cK_{\gb ,\bdf}}\|_{\bbL_1} \leq \gd ) 
\geq \frac14\, \jpm(\cA\cap{\tilde \cA})
= \frac14 \, \frac{\RCbdf{\wired}{\abs\bdf}(\cA\cap{\tilde \cA})}
{\RCbdf{\wired}{\abs\bdf}(\frI)}\,.
\label{eq_goal}
\end{align}
The denominator can be easily estimated by \eqref{1_Delta_beta_eta_fk}
\begin{equation}
\label{eq boundary surface tension}
\RCbdf{\wired}{\abs\bdf}(\frI) =  
\exp \big(- \Delta_{\gb ,|\bdf |} \, N^{d-1} + \smallo(N^{d-1}) \big) \, .
\end{equation}
On the other hand, as we are trying to bound probabilities from
below, even under the conditioning on $\tilde{\cA}$  
one can still apply FKG arguments (non-crossing is a decreasing event)
and use  
 \eqref{1_tau_beta_fk} and  \eqref{1_Delta_beta_eta_fk} to decouple between 
different $\hat B_i$-events which constitute $\cA$ in \eqref{calA}. 
Therefore, since, 
$\RCbdf{\wired}{\abs\bdf}({\tilde \cA})\geq {\rm e}^{-c_3\delta N^{d-1}}$, we are 
entitled to conclude that 
\begin{equation}
\label{eq lower bound surface tension}
\RCbdf{\wired}{\abs\bdf}(\cA\cap{\tilde \cA}) \geq
\exp \left\{ - N^{d-1}\big(
\sum_{\vec{n}_i\neq -\vec{e}_d}\tau_\gb (\vec{n}_i) |\widehat B^i |
\, +\,
\Delta_{\gb ,|\bdf |} \sum_{\vec{n}_i= -\vec{e}_d} |\widehat B^i |
\big)\,  
- \,\big(c_3\gd +\smallo (\gd)\big) N^{d-1} \right\} .
\end{equation}
%
Combining  
\eqref{eq lower bound surface tension} 
 with the estimates 
\eqref{eq_goal} and \eqref{eq boundary surface tension};
the approximation property \eqref{W_pm_approx}; the 
choice of the polyhedral set $\cK_{\gb ,\bdf}^\gd (v(m))$ through 
Lemma~\ref{thm approx} and the representation \eqref{functional_g} of 
$\cW_{\gb ,\bdf}^{\pm}$, we arrive  to 
the claim of Proposition~\ref{prop 2}.
\end{proof}

\subsection{Approximation result II}
Approximation results used in the proof of 
 the upper bound should apply to any 
$v\in \BV$. 
On the other hand, on the level of the precision provided by the 
${\Bbb L}_1$-theory 
the underlying probabilistic estimates
are less delicate  than those needed for the  
 proof of  the lower bound for the distinguished 
almost optimal polyhedral shapes of the preceding subsection. 
In fact, global bulk relaxation properties play no r\^{o}le at all,
 and the ${\Bbb L}_1$-upper bound is just a coarse estimate 
on the surface tension
 price of the localized interface. In particular, the approximation
 construction of \cite{Bo} suffices and the proof of the corresponding result
 in the latter paper literally goes through. 
As in the case of the lower bound let us concentrate on the 
more difficult case of the negative boundary magnetic field $\bdf\leq 0$.
 We claim:  
\begin{thm}
\label{theo ABCP} 
Let $v\in \BV$. For any $\gd$ positive, there exists $h$ positive and a finite
number 
  $\ell$ of disjoint parallelepipeds ${\widehat R^1}, \dots, \widehat
R^{\ell}$   with basis $\widehat B^1, \dots, \widehat B^\ell$ of  size 
$h$ and 
height $\gd h$ such that:\\
{\bf II\,a)}\ For every $i=1,...,\ell$ 
 the base $\widehat B^i$ divides $\widehat R^i$ into two 
parallelepipeds  $\widehat R^{i,+}$ and $\widehat R^{i,-}$ and we denote by
$\vec{n}_i$ the normal to $\widehat B^i$. 
Either $\widehat B^i \cap \partial \uBox = \emptyset$ and $\widehat
R^i$ is included in $\uBox$, or $\widehat B^i$ is
included in $\partial \uBox \setminus \cP_d$ and only $\widehat
R^{i,-}$ is included in $\uBox$, or $\widehat B^i$ is included in $\cP_d$ 
and only $\widehat R^{i,+}$ is included in $\uBox$\\
\noindent
{\bf II\,b)}\  The parallelepipeds ${\widehat R^1}, \dots, \widehat
R^{\ell}$ approximate the reduced boundary $\partial^*_{g^{\pm}}v$ (see 
Subsection~\ref{ssec_BV}) in the following sense:
$$
\int_{\widehat R^i} | \cX_{\widehat R^i}(x)  - v\vee g^\pm (x)| \, dx \leq \gd \, 
\vol(\widehat R^i),
$$
where
 $\cX_{\widehat R^i} = 1_{\widehat R^{i,+}} - 1_{\widehat R^{i,-}}$ 
and the volume of $\widehat R^i$ is ${\vol(\widehat R^i)} = \gd h^d$. Furthermore, 
$$
\Big| \sum_{\widehat B^i\cap\cP_d = \emptyset }\tau_\gb (\vec{n}_i)| \widehat B^i |
\, +\,
\Delta_{\gb ,\abs\bdf}\sum_{\widehat B^i\subseteq \cP_d } |\widehat B^i| 
\, 
 - \,\cW_{\gb,\abs\bdf}^{\pm}(v) \Big| \leq \gd,
$$
\end{thm}

As we have mentioned, 
Theorem~\ref{theo ABCP} (see figure \ref{fig_approx}) is proved as in
the standard ($+$ b.c.) case, see~\cite{Bo}.
\subsection{The upper bound}

\begin{pro}
\label{prop 3}
We fix $\gb \in \frB$ and $\bdf < 0$. For all $v$ in $\BV$ such that 
$\cW_{\gb,\bdf}(v)$ is finite, one can choose $\gd_0 = \gd_0 (v)$,
such that uniformly in $\gd < \gd_0$
\begin{eqnarray*}
\limsup_{N \to \infty} \; \frac{1}{N^{d-1}}
\max_{K_0(\gd)  \leq K \leq N^\nu}
\log \Isbdf{+}{N}{\bdf}  
\left( \| \cM_{N,K} - m^* v \|_{\bbL_1} \leq \gd \right) \leq 
- \intSTbdf{\bdf}(v) 
+ o(\gd) \, ,
\end{eqnarray*}
where the function $o(\cdot)$  depends only on $\gb$ and $v$ and vanishes 
as $\gd$ goes to 0.
\end{pro}

We fix $K$ large enough and $\gd > 0$.
Using the approximation Theorem, there is $\gd_0 < \gd$ such that for
any $\gd ' < \gd_0$
\begin{align*}
\limsup_{N\ra\infty} \frac1{N^{d-1}} \log&
\Isbd{+}{N} \left( 
\| \frac{1}{m^*} \cM_{N,K} - v \|_{\bbL_1} \leq \gd ' \right)\\
&\leq \limsup_{N\ra\infty} \frac1{N^{d-1}} \log
\Isbdf{\pm}{N}{\abs\bdf} \left(
\frac{\cM_{N,K}}{m^*} \in \inter_{i=1}^\ell
\cV(\widehat R^i,2\gd\vol \widehat R^i) 
\right)\\
&= \limsup_{N\ra\infty} \frac1{N^{d-1}} \log \jpm
\left( \frac {\cM_{N,K}}{m^*} 
\in \inter_{i=1}^\ell \cV(\widehat R^i,2\gd\vol \widehat R^i) \right),
\end{align*}
where $\cV( \widehat R^i , \gep)$ is the $\gep$-neighborhood of 
$\cX_{\widehat R^i}$ 
\begin{eqnarray*}
\cV (\widehat R^i, \gep) = \left\{ v^\prime \in \bbL_1( \uBox )
\ \big| \quad  \int_{\widehat R^i} | v^\prime(x) - \cX_{\widehat R^i}(x) |  
\, dx 
\leq \gep  \right\}.
\end{eqnarray*}

According to Lemma \ref{lem coarse graining}, there exist $K_0$ and
$\gz$ such that
\begin{align*}
\limsup_{N\ra\infty} \frac1{N^{d-1}} & \max_{K_0 \leq K \leq N^\nu} 
\log \jpm
\left(
\frac {\cM_{N,K}}{m^*} \in \inter_{i=1}^\ell \cV(\widehat R^i,2\gd\vol \widehat
R^i)
\right) \\
& \leq \limsup_{N\ra\infty} \frac1{N^{d-1}} 
\max_{K_0 \leq K \leq N^\nu}
\log \jpm
\left(
u^\gz_{N,K} \in
\inter_{i=1}^\ell \cV(\widehat R^i,3\gd\vol \widehat R^i)
\right) .
\end{align*}

The first step to extract the surface tension factor is to localize in each
$R^i_N$ a mesoscopic interface
(where $R^i_N$ is the microscopic counterpart of $\widehat R^i$).
This is done by means of a surgical procedure :
in each $R^i_N$, the $\bbL_1$-constraint  $u^\gz_{N,K} \in 
\cV(\widehat R^i , 3 \gd \vol(\widehat R^i))$ enables to find two sections
in $R^{i,+}_N$ and in $R^{i,-}_N$  (the minimal sections) containing only 
a small portion of the mesoscopic interface.
So that up to a small cost one can rearrange the configurations in
these sections in order to identify clearly the location of the interface.
At that point the microscopic structure of Pisztora's coarse graining
becomes effective (see Remark \ref{0-block}) and the presence of a 
mesoscopic interface is enough to ensure that the sets 
$\bnd^{\rm top}R^{i'}_N$ and $\bnd^{\rm bot} R^{i'}_N$ are disconnected 
on the microscopic level, where ${\widehat R^i} \, '$ is the parallelepiped 
included in $\widehat R^i$ with basis $\widehat B^i$ and height 
$\frac{\gd}{2} h$. Therefore, introducing the set 
\begin{eqnarray*}
\cA = \left\{ \go \in \frI \ : \ 
\exists \gs \ {\rm such \  that}  \  u^\gz_{N,K} (\gs,\go) \in 
\bigcap_{i =1}^{\ell} \,  \cV(\widehat R^i , 3 \gd \vol(\widehat R^i))
\right\} \, ,
\end{eqnarray*}
we obtain
\begin{equation}
\label{cA}
\jpm(\cA) = \RCbdpmf{\abs\bdf}(\cA) \leq e^{C\gd N^{d-1}}\,
\RCbdpmf{\abs\bdf} \left( \inter_i\{\bnd^{\rm top}R^{i'}_N \nlra \bnd^{\rm
bot} R^{i'}_N\} \right) \,,
\end{equation}
The constant $C$ in the error term is proportional to the perimeter
of $v$.
It is important to note that the surgical procedure is insensitive to the
boundary effects and applies equally well near the wall $\cP_d$ or in the bulk.

Equation \eqref{cA} can then be bounded by
\begin{equation}
\RCbdpmf{\abs\bdf} \left( \inter_i\{\bnd^{\rm top}R^{i'}_N \nlra \bnd^{\rm bot}
R^{i'}_N\} \right)
\leq \frac{\RCbdf{\wired}{\abs\bdf}(\inter_i\{\bnd^{\rm top}R^{i'}_N \nlra
\bnd^{\rm bot} R^{i'}_N\})} {\RCbdf{\wired}{\abs\bdf}(\frI)}\,.
\end{equation}
As in the proof of the lower bound, the denominator can be estimated 
by \eqref{1_Delta_beta_eta_fk}.
In order to recover the surface tension in each box $R^i_N$, we fix the
boundary conditions outside each box $R^i_N$.
At that point the relaxed definition of surface tension is crucial
and we  apply  Lemma
\ref{lem surface tension} or \ref{lem wall free energy} depending whether
the box $R^i_N$ lies on the wall or in the bulk.
This provides the following upper bound
\begin{equation}
\RCbdf{\wired}{\abs\bdf} ( \cA ) \leq 
\exp \left\{  - N^{d-1} 
\left(
\sum_{\widehat B^i \cap\cP_d = \emptyset} 
\tau_\gb (\vec{n}_i) \, |\widehat B^i|\,
+\, \Delta_{\gb,\abs\bdf}
\sum_{\widehat B^i \subseteq \cP_d} |\widehat B^i|
\right) \,
+\,  C \gd N^{d-1} 
\right\} \, .
\end{equation}
Using the approximation Theorem \ref{theo ABCP}, we, thereby, obtain:
\begin{equation}
\RCbdf{\wired}{\abs\bdf} ( \cA ) \leq \exp \big\{ - N^{d-1} \cW_{\gb,\abs\bdf}^{\pm}
(v)  + \smallo(\gd) N^{d-1} \big\} \, .
\end{equation}
In view of \eqref{eq boundary surface tension} and the expression 
 \eqref{functional_g} for the functional $\cW_{\gb,\abs\bdf}^{\pm}$  the proof of  
Proposition~\ref{prop 3} is completed.

\qed

\subsection{Proof of Theorems~\ref{theorem_A} and \ref{theorem_B}}

The Theorems~\ref{theorem_A} and \ref{theorem_B} follow from
Propositions \ref{prop 1}, \ref{prop 2}, \ref{prop 3}. 
We will focus on the derivation of Theorem~\ref{theorem_B} in the case of
negative boundary field, the other cases can be deduced easily.\\

The first step is to derive the following lower bound for $m \in ]
\bar m (\gb,\bdf), m^*_\gb[$
\begin{eqnarray}
\label{lower bound eq}
\liminf_{N \to \infty} \frac{1}{N^{d-1}}
\log \Isbd{+}{\BoxxN}\left( {\bf M}_{\BoxxN}\leq m\right)
\geq - w^\star_{\gb ,\bdf} (m) \, , 
\end{eqnarray}
where $w^\star_{\gb ,\bdf} (m)$ was introduced in \eqref{claim_A}. 

Fix $\gep$ in $]0,m^*_\gb -  \bar m(\gb,\bdf)[$.
For any $\gd '>0$ small enough and $\ga$ in $(0,1/d)$, one 
has
\begin{eqnarray*}
\Isbd{+}{\BoxxN} \left( {\bf M}_{\BoxxN}\leq m\right) \geq 
\Isbd{+}{\BoxxN} \left( \| \cM_{N,N^\ga} - m^* 
\1_{\cK_{\gb,\bdf}(v(m - \gep ))} \|_{\bbL_1}
\leq \gd ' \right) \, .
\end{eqnarray*}
By using Proposition \ref{prop 2} and letting $\gd '$ go to 0, we see 
that 
\begin{eqnarray*}
\liminf_{N \to \infty} \frac{1}{N^{d-1}}
\log \Isbd{+}{\BoxxN} \left( {\bf M}_{\BoxxN}\leq m\right) \geq 
- w^*_{\gb,\bdf} (m - \gep)   \, .
\end{eqnarray*}
As the minimizers $\cK_{\gb ,\bdf} (v(m))$ are obtained by
dilatation of the same unnormalized Winterbottom shape 
$\cK_{\gb ,\bdf}$, it follows that
\begin{eqnarray*}
\label{continuity intSTbdf}
\lim_{\gep \to 0} \,w^*_{\gb,\bdf} (m - \gep)\,=\, 
w^*_{\gb,\bdf} (m) \, .
\end{eqnarray*}
This concludes the proof of \eqref{lower bound eq}.


\medskip

In order to derive the upper bound, we define for a given
$\gd >0$
\begin{eqnarray*}
F_\gd = \left\{ u \in \bbL^1(\uBox) \; \Big| \;
\inf_{x} 
\|u - \1_{x +\cK_{\gb,\bdf}(v(m))}\|_{\bbL_1}\geq \gd;
\quad  \int_{\uBox} u (y) \, {\rm d}y \leq \frac{m}{m^*} 
\right\} \, .
\end{eqnarray*}
For a given $a>0, \gd' >0$,
\begin{eqnarray*}
\Isbd{+}{\BoxxN} \left( \frac{1}{m^*} \cM_{N,K} \in F_\gd \right) 
& \leq
\Isbd{+}{\BoxxN} \left( \frac{1}{m^*} \cM_{N,K} \in F_\gd 
\cap \cV\left(\cC_{g^\pm}(a),\gd' \right) \right)\\
& \qquad +
\Isbd{+}{\BoxxN} \left( \frac{1}{m^*} \cM_{N,K} \not \in
\cV\left(\cC_{g^\pm}(a),\gd' \right) \right)
\, .
\end{eqnarray*}

According to Proposition \ref{prop 1}, there is $K_0$ depending on
$\gd'$ such that uniformly in $a$
\begin{eqnarray*}
\limsup_{N \to \infty} \;   \frac{1}{N^{d-1}}
\max_{K_0 \leq K \leq N^\nu}
\log 
\jpm 
\left( 
\frac{1}{m^*} \cM_{N,K}  \not \in \cV\left(\cC_{g^\pm}(a),\gd '\right)
\right) 
\leq - \frac{C_\gb}{K_0^{d-1}} a + \left|\taubd\right| \, .
\end{eqnarray*}

We fix $\gep >0$ and $a$ much larger than $\frac{K_0^{d-1}}{C_\gb} 
(w^*_{\gb,\bdf} (m) + \left|\taubd\right|)$.
The exponential tightness will enable us to consider only a finite number
of macroscopic configurations:

Since $\cC_{g_{\pm}}(a)\cap F_\gd$ is compact in ${\Bbb L}_1 (\uBox )$, we, using 
the upper bound Proposition~\ref{prop 3}, can cover it (and hence 
also $\cV\left(\cC_{g^\pm}(a)\cap F_\gd,\gd '\right)$ for $\gd'$ sufficiently small)
 with a finite net of 
neighborhoods $\big( \cV(u_i,\gep_i) \big)_{i \leq \ell}$ and choose a finite scale
 $K_1 (\gep)$ such that for every $i \leq \ell$ 
\begin{eqnarray*}
\limsup_{N \to \infty} \; \frac{1}{N^{d-1}}
\max_{K_1 (\gep)  \leq K \leq N^\nu}
\log \Isbdf{+}{N}{\bdf}  
\left( 
\frac{1}{m^*} \cM_{N,K}  \in \cV (u_i,\gep_i)
\right) 
\leq  - \intSTbdf{\bdf}(u_i) + \gep \, .
\end{eqnarray*}
%
By the stability property \eqref{stability}, 
\begin{eqnarray*}
\inf_{u \in F_\gd}  \intSTbdf{\bdf}(u) > w^*_{\gb ,\bdf}(m) \, .
\end{eqnarray*}
Therefore, by choosing $\gep$ small enough, 
 we can find
a mesoscopic scale $K_2$  
and a positive constant $c = c(\delta )>0$,
 such that
\begin{eqnarray*}
\limsup_{N \to \infty} \;   \frac{1}{N^{d-1}}
\max_{K_2 \leq K \leq N^\nu}
\log \Isbd{+}{\BoxxN} \left( \frac{1}{m^*} \cM_{N,K} \in F_\gd \right) 
\leq - \, w^*_{\gb,\bdf}(m) -c(\gd ) .
\end{eqnarray*}
The above inequality combined with the lower bound \eqref{lower bound eq}
leads to the result. \qed

\medskip

\noindent
{\bf Acknowledgments:} 
TB would like to thank A. Chambolle for many interesting discussions and
suggestions about Geometric measure Theory. 
TB and YV also acknowledge the warm hospitality of Technion, where
part of this work was done.

\end{document}